\documentclass[12pt]{article}
\usepackage{amsmath,amssymb,amsthm}
\usepackage{psfig}
\usepackage{subfigure}

\newtheorem{thm}{Theorem}%[section] (If you want theorem numbered
\newtheorem{lemma}{Lemma}%               with section number.  Same
\newtheorem{cor}{Corollary}%       goes for lemmas, etc.)
\newtheorem{prop}{Proposition} %--> \begin\end{theorem,lemma,...}
\newtheorem{rem}{Remark}

\newtheorem{defin}{Definition}

\newcommand{\x}{\times}

\newcommand{\Ga}{\Gamma}

\newcommand{\ga}{\gamma}
\newcommand{\ep}{\epsilon}

\newcommand{\C}{\mathbb C}

\newcommand{\e}{\varepsilon}

\newcommand{\mor}{morphism}

\newcommand{\p}{\partial}
\newcommand{\Dt}{D_t}
\newcommand{\var}{\varphi}
\newcommand{\z}{\mathbb Z}
\newcommand{\pl}{{\mathbb R}^2}

\begin{document}
\setlength{\baselineskip}{16pt}

\title{Invariants of tangles with flat connections in their complements.I.
Invariants and holonomy $R$-matrices}

\author{R. Kashaev,  N.Reshetikhin}
\maketitle

\begin{abstract}The notion of holonomy $R$-matrices
is introduced. It is shown how to define invariants of tangles
with flat connections in a principle $G$-bundle of the complement
of a tangle using holonomy $R$-matrices. 
\end{abstract}

\section{Introduction}

A breakthrough in the theory of invariants of knots
came with the discovery of the Jones
polynomial \cite{J} and its generalizations (HOMFLY and Kauffman).
The construction of the Jones invariant based on a solution to the
Yang-Baxter equation was given in \cite{J-1} and was generalized
to HOMFLY polynomial and to more general case in \cite{T}. It involved
a solution to the Yang-Baxter equation which satisfies certain extra
condition.  Then it was shown
in \cite{R} \cite{RT} that such invariants can be obtained from quantized
universal algebras for any simple Lie algebra.
In \cite{R-1} it was proven that the non-degeneracy
and the cross-nondegeneracy of a solution of the Yang-Baxter equation
is enough to construct the invariants.

In this paper we re-examine the construction of invariants of
links from \cite{RT} and \cite{R-1} (for details see \cite{T-1}.
Using systems of holonomy
$R$-matrices we construct the invariant of tangles in $\pl\times
I$ ( here $I$ is a unit interval) with a flat connection in a
principle $G$-bundle ( where $G$ is a simple complex algebraic
group) over the complement of a tangle. We show that this
invariant depends only on a gauge class of this connection.

In the follow-up paper we will show that quantized universal
$R$-matrices provide examples of holonomy $R$-matrices.

In the second section of this paper we describe the moduli space
of flat connections in the complement to a tangle in ${\mathbb
R}^3$. In the third section we describe the category of tangles
with flat connections in the complement. This construction is very similar
to $\pi$-tangles introduced and studied in \cite{T-2}. Then in section 4
we construct invariants of tangles with flat connections
in the complement using non-degenerate solutions to
the holonomy Yang-Baxter equations. In section 5
we show that these invariants are gauge invariant, and in particularly
for links they can are functions (sections of a line bundle)
over the moduli space of flat connections in the complement
of tangle. Some proofs are moved to the Appedix.

This paper is one of a series of papers in which we will analyze
invariants of tangles and 3-manifolds related to quantized universal
enveloping algebras with large center. Some outline of this program is
given in the conclusion (section 6).

The main results of this paper and of the part II were announced at the conference
"Graphs and Patterns in Mathematics and Physics". Both authors were
partly supported by the NSF grant DMS-0070931 and by the CRDF grant RM1-2244.

\section{Flat connections in the complement of a tangle}

\subsection{The fundamental group of the complement of a tangle}\label{pigr}

Let $I\equiv[0,1]$ be the unit closed interval.
A geometric tangle $t\subset {\pl}\x I$ is the image of an
embedding
\[
\underbrace{I\sqcup\ldots\sqcup I}_{k\ \mathrm{ times}}
\sqcup \underbrace{S^1\sqcup\ldots\sqcup S^1}_{l\ \mathrm{ times}}
\to \pl\x I
\]
with oriented components such that
\[
\p t\subset {\pl}\x\p I={\pl}\x\{0,1\}
\]
and $\p t\cap {\pl}\x \p I= t\cap {\pl}\x\p I$.

Let $(x,y,z)$, $x,y\in{\mathbb
R}$, $z\in I$, be Cartesian coordinates on ${\pl}\x I$.
Denote by $\Dt$ the image of $t$ under the projection
\begin{equation}\label{eq:proj}
p\colon {\pl}\x I\ni(x,y,z)\mapsto(y,z)\in\pl
\end{equation}
We call geometric tangle $t$ {\it restricted}
if it is based on lines
$L_\pm=\{(0,y,(1\pm1)/2)\}$, i.e. $\p t\in L_+\cup L_-$.
Denote by
\[
\p_{\pm}t\equiv \p t\cap L_\pm,\quad
\p_{\pm}\Dt\equiv p(\p_\pm t)
\]
the boundary components of $t$ and $\Dt$ which belong to $L_{\pm}$ and their
images in $\Dt$,
respectively.

\begin{figure}
\centering
\mbox{\subfigure[Positive intersection]{\psfig{figure=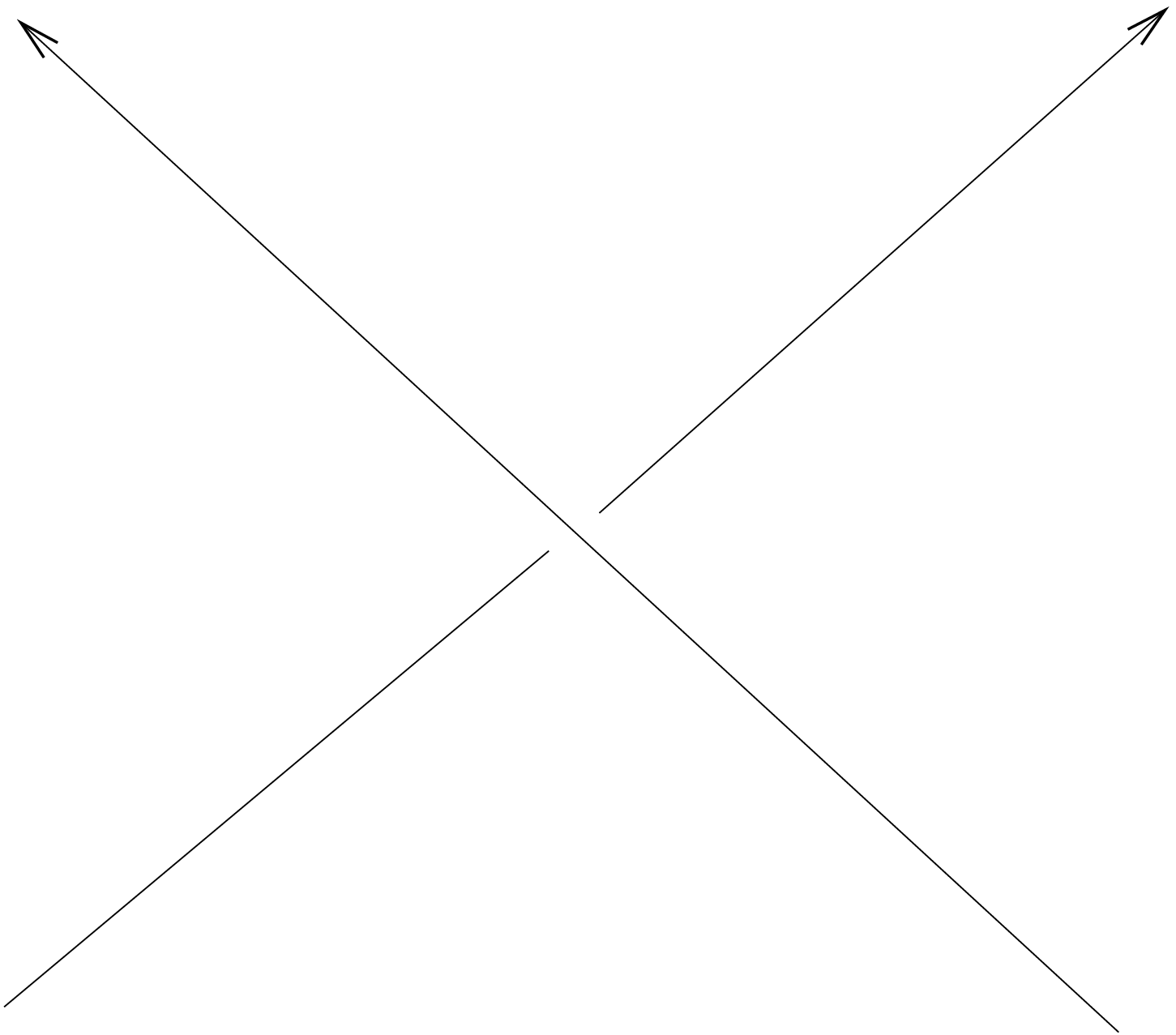,width=1.2in}}\qquad
      \subfigure[Negative intersection]{\psfig{figure=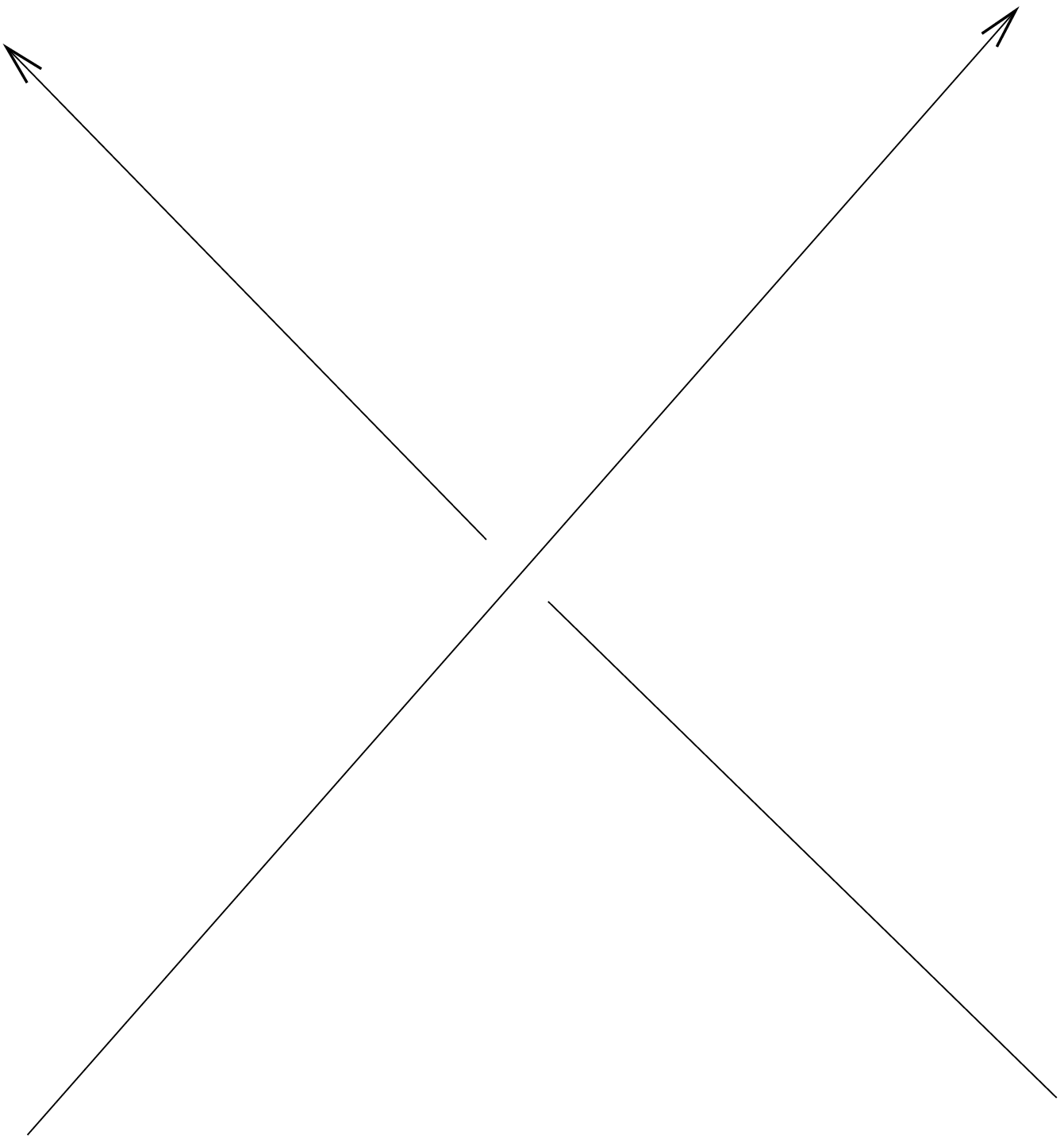,width=1.2in}}}
\caption{\label{pos-neg} }
\end{figure}

Projection $\Dt$ is called {\it regular} if it has only double points at which
segments intersect transversely, all critical points in
$z$-direction are
non-degenerate, and non of the critical points coincide with any
of the double points. A regular projection $\Dt$ will be called
\emph{diagram} of $t$. The vertices of a diagram are its double points,
critical points in $z$-direction and boundary points, while edges are the
segments running between the vertices.
A vertex of $D_t$ is called  \emph{positive or $(+)$ crossing} if it is
a double point and the angle between lower
and upper components is positive.
If the angle is negative, the vertex is called \emph{negative or $(-)$
crossing}, see Fig.~\ref{pos-neg}.

\begin{figure}
\centerline{\psfig{figure=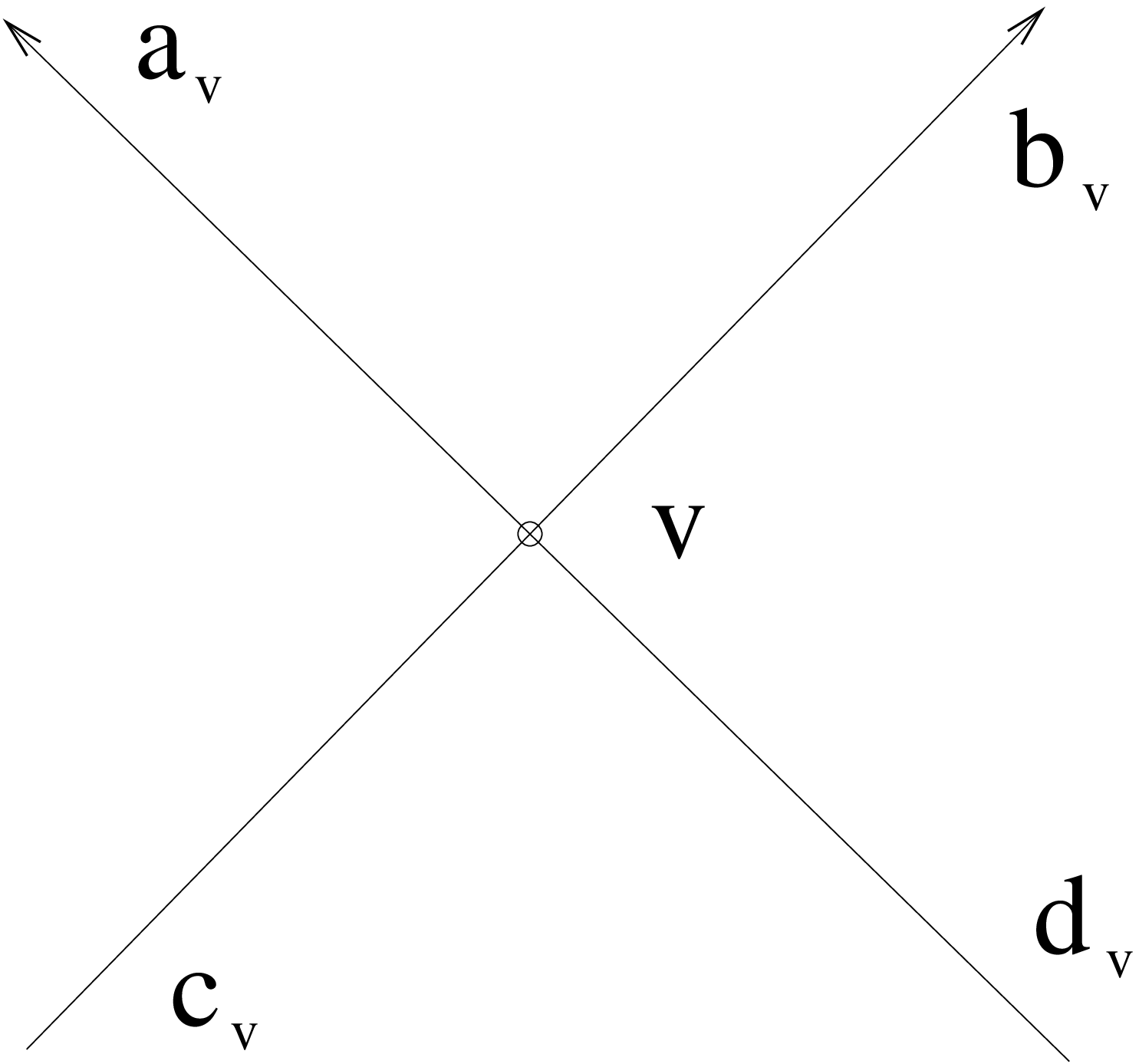,width=2.0in}}
\caption{\label{names}}
\end{figure}

\begin{figure}
\centerline{\psfig{figure=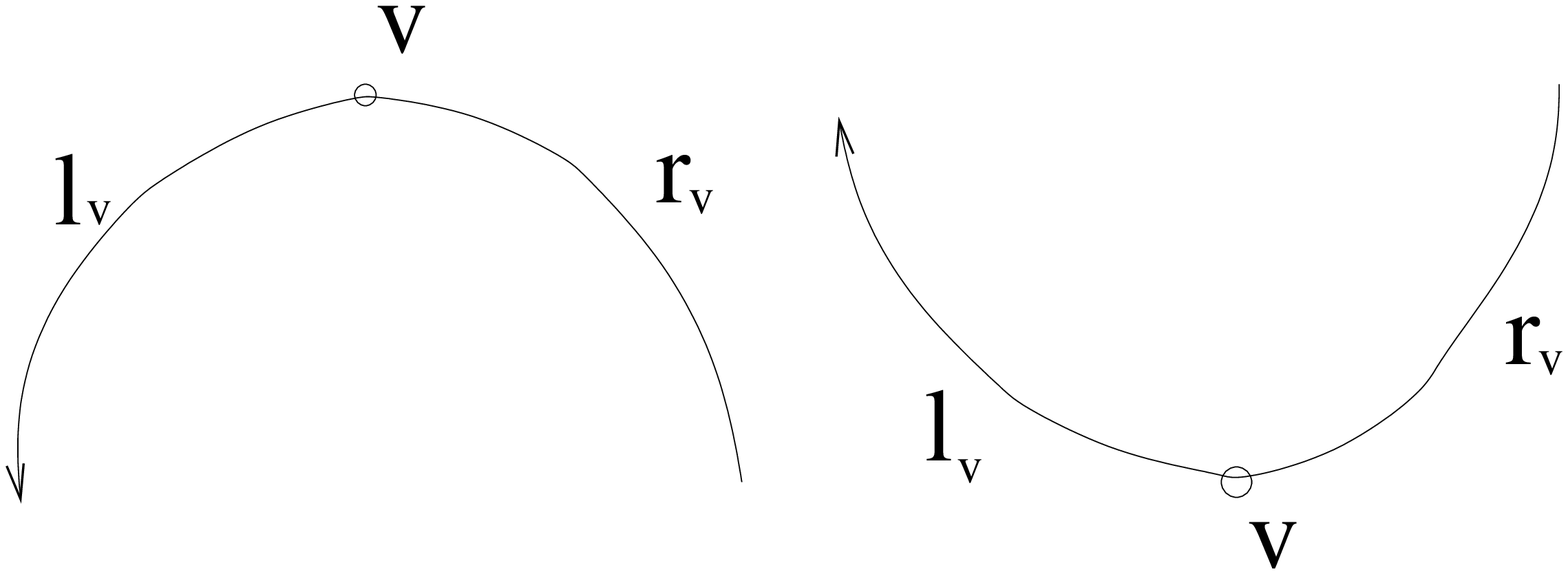,width=2.0in}}
\caption{\label{eval-inj}}
\end{figure}

\begin{figure}
\centerline{\psfig{figure=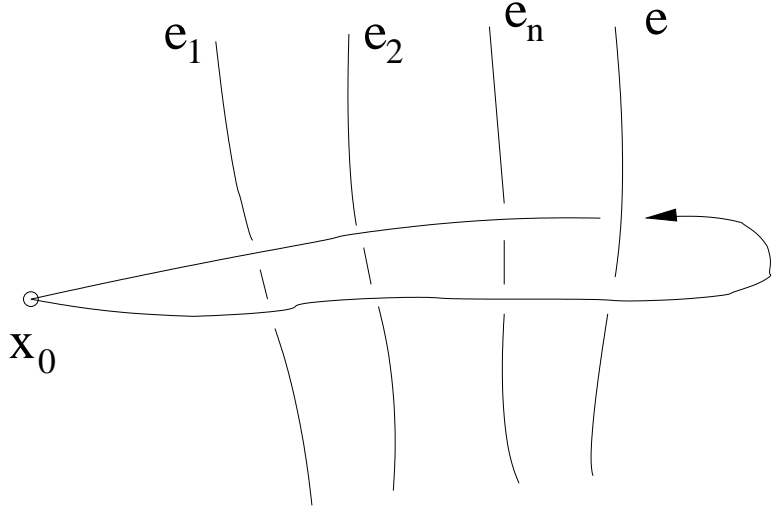,width=4.0in}}
%\centerline{\psfig{figure=fig-mod.eps,width=4.0in}}
\caption{\label{mon} The choice of the monodromy path }
\end{figure}

For each double point $v$ of diagram $D_t$ denote by
$a_v,b_v,c_v,d_v$ the adjacent edges in accordance with Fig.~\ref{names},
and for each critical point $v$
denote by $l_v, r_v$ the adjacent edges as in Fig.~\ref{eval-inj}.
%\begin{figure}
%\centering
%\mbox{\subfigure[]{\psfig{figure=fig-3a.eps,width=1.2in}}\qquad
%      \subfigure[]{\psfig{figure=fig-3b.eps,width=1.2in}}}
%%\centerline{\psfig{figure=fig-1.eps,height=1.2in,width=4.0in}}
%%\caption{\label{pos} Positive intersections}
%\end{figure}
Define group $\Pi(D_t)$ generated by edges of $D_t$ with the
following defining relations:
\begin{itemize}
\item for each double point $v$
\[
\begin{array}{rl}
 (+)\mbox{ crossing:}&a_v=d_v,\ b_v=d_v^{-1}c_vd_v \\
(-)\mbox{ crossing:}&b_v=c_v,\ a_v=c_vd_vc_v^{-1}
\end{array}
\]
\item for each critical point $v$: $l_vr_v=1$.
\end{itemize}
The isomorphism class of this group does not
depend on the diagram representing $t$ and is
isomorphic to the fundamental group of the complement to $t$:
\[
\Pi(D_t)\simeq \pi_1(\pl\times I\backslash t)
\]
In fact, the definition of  $\Pi(D_t)$ concides with the
Wirtinger presentation of  $\pi_1(\pl\times I\backslash t)$.
One can fix an isomorphism by choosing a base point $x_0$ in
$\pl\times I\backslash t$ and associating to each edge $e$ in $D_t$
element $\Ga^{x_0}_e$ of $\pi_1(\pl\times I\backslash t,x_0)$ which is
represented by the loop at $x_0$ with projection as in
 Fig.~\ref{mon}. We shall call such correspondence
a \emph{monodromy} map.

Associated with boundary of $D_t$, group $\Pi(\Dt)$ has two natural subgroups
$\Pi(\p_{\pm}\Dt)$ generated by edges intersecting with corresponding
boundary components. One has the following commutative diagram:
\[
\begin{array}{ccccc}
\Pi(\p_- \Dt) & \hookrightarrow& \Pi(\Dt)&\hookleftarrow& \Pi(\p_+\Dt)\\
\uparrow&&\wr|&& \uparrow\\
\pi_1({\pl}\x\{0\}\setminus\p_- t,x^-_0)&\rightarrow&
\pi_1({\pl}\x I\setminus t,x_0)&\leftarrow&\pi_1({\pl}\x \{1\}\setminus
\p_+ t,x^+_0)
\end{array}
\]
where the hookarrows are embeddings while
the arrows are natural group homomorphisms, and we denote
$x^\pm_0=(x,y,(1\pm1)/2)$ the projections of the base point
$x_0=(x,y,z)$. The base point is chosen such that it projects
to the plane of the diagram $D_t$ to the far left of it
(in the $y$-direction).

\subsection{The moduli space of flat connections in the complement to a
tangle}\label{modspace}

For a  Lie group $G$ let $A$ be the
1-form corresponding to a flat connection in a principle $G$-bundle over the
complement to tangle $t$ in ${\mathbb R}^2\times I$. Such form
defines a group homomorphism $\pi_1({\mathbb R}^2\times I\backslash
t, x_0)\to G$ which maps the element generated by path $\Ga$
to the holonomy of $A$ along $\Ga$. Because we have the isomorphism
$\pi_1({\mathbb R}^2\times I\backslash t,x_0)\simeq\Pi(D_t)$ a
flat $G$-connection also defines a group homo\mor
\[
\gamma\colon \Pi(D_t)\ni e\mapsto \gamma_e\in G
\]
We denote the space of all such homo\mor s  as $\Pi(D_t, G)$.
The group $G$ acts on $\Pi(D_t, G)$ by conjugations
\[
g: \ga_e\mapsto g\ga_e g^{-1}
\]
In terms of $\pi_1({\mathbb R}^2\times I\backslash t,x_0)\simeq\Pi(D_t)$
this action corresponds to the changes of the trivialization of the
$G$-bundle at the base point, i.e. to gauge transformations.

Thus, the gauge class of the flat connection $A$ defines an
element in the coset space $\Pi(D_t, G)/G$. This gives a
combinatorial description of the space of isomorphism classes
of representations of the fundamental group of
${\mathbb R}^2\times I\backslash t$ in $G$ ( or equivalently,
of the modulai space of flat connections):
\[
\Pi(D_t, G)/G\simeq \{\pi_1({\mathbb R}^2\times I\backslash
t, x_0)\to G\}/G
\]
Notice that this description of the isomorphism classes of representations
of the fundamental group works also in cases when $G$
is not a Lie group.

Notice that the homomorphism $f: \pi_1(\pl\times I\backslash, x_0)\to G$
does not define a flat connection in the complement of the tangle.
It only defines a gauge class of a flat connection.

\begin{figure}
\centerline{\psfig{figure=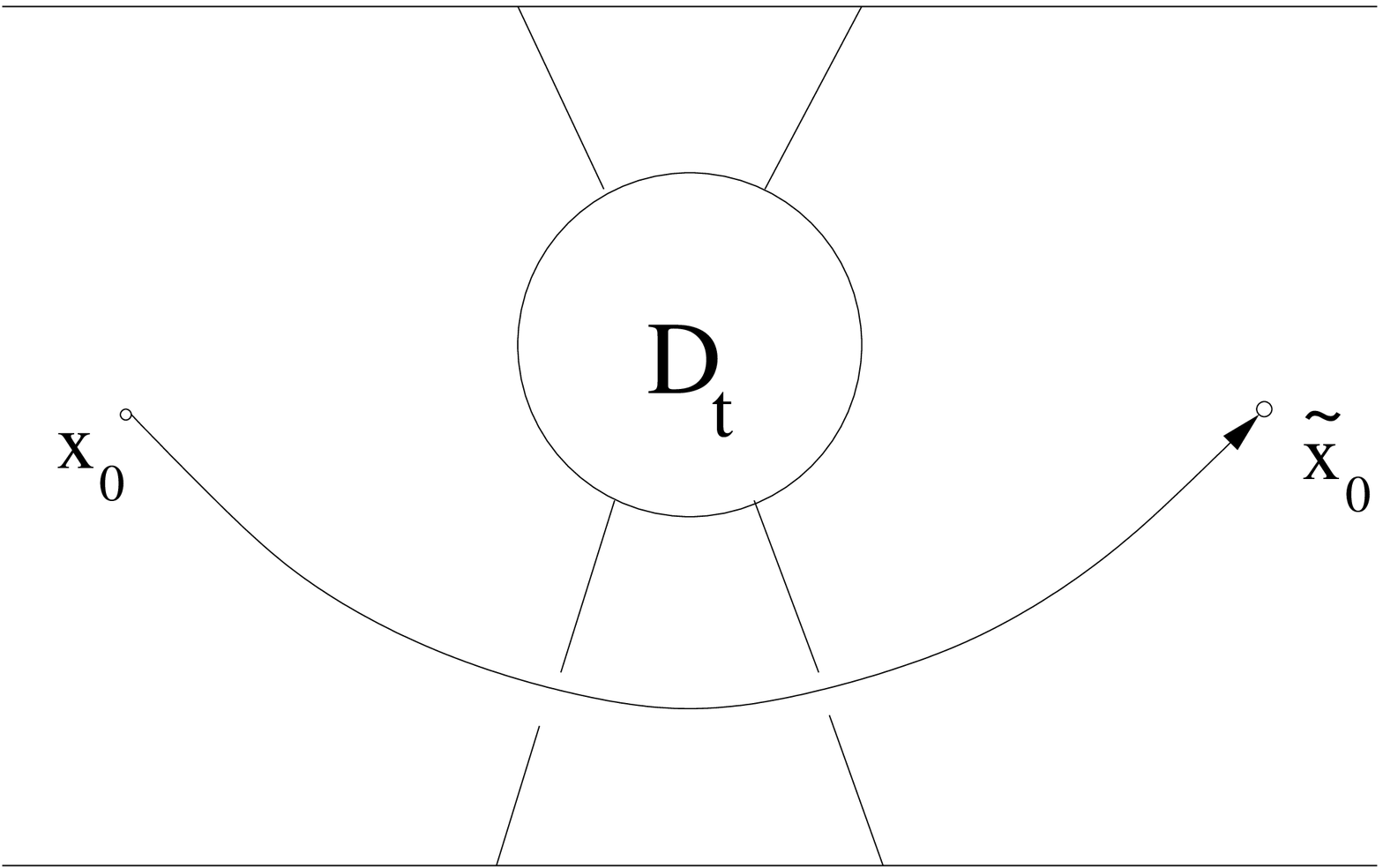,width=4.0in}}
%\centerline{\psfig{figure=fig-mod.eps,width=4.0in}}
\caption{\label{fig-ramif}}
\end{figure}

For what will follow we need slightly
more then a representation of the fundamental group, but less then
a flat connection. Let $\tilde{x}_0$ be a point which projects to
the far right from $D_t$.
A pair $(f,g)$ where $f$ is a representation of the
fundamental group and $g\in G$ is called a {\it ramified representation
of the fundamental group} if $g$ has a meaning of the holonomy of
a flat connection whose gauge class is defined by $f$ along the path
connection $x_0$ and $\tilde{x}_0$ in such a way that it projects "over"
the diagram $D_t$ (see Fig. \ref{fig-ramif}).

Let $\{g_e\}_{e\in E(D_t}, g_r$ be a ramified representation of the
fundamental group.
The gauge group acts on it as
$g_e\mapsto gg_e g^{-1}, \ g_r\mapsto gg_r h^{-1}$. Thus the action
of the gauge group on ramified representations reduces to the action of
$G\times G$. The isomorphism class does not depend on the ramification.

\subsection{Factorizable groups and Lie groups}

We say that group $G$ is factorizable into two
subgroups $G_\pm\subset G$ if any element $g\in G$
can be represented in a unique way as
\begin{equation}\label{la-factor}
g=g_+g_-^{-1}=\bar{g}_-^{-1}\bar{g}_+
\end{equation}
where $g_\pm, \bar{g}_\pm \in G_\pm$.

Let $G$ be a Lie group with Lie algebra $\mathfrak{g}$. We
say that Lie algebra $\mathfrak{g}$ is factorizable into
two Lie subalgebras $\mathfrak{g}_\pm$ if for each $x\in
\mathfrak{g}$ there exists unique pair $x_\pm\in \mathfrak{g}_\pm$
such that
\[
x=x_+ -x_-
\]
A Lie group with factorizable Lie algebra is almost factorizable
(i.e. factorizable on an open dense subset). It will be called
factorizable Lie group. We will be especially interested in the
case when $G$ is an algebraic Lie group with factorizable Lie
algebra. Then it is factorizable on a Zariski open subset of $G$.

Another important case is the conditional factorization. In
conditionally factorizable Lie algebras factorization
(\ref{la-factor}) is unique only under extra conditions on $x_\pm$.
For example, fixing a Borel subalgebra ${\mathfrak
b}_+$ in a simple Lie algebra $\mathfrak g$ makes $\mathfrak g$
conditionally factorizable with ${\mathfrak g}_\pm={\mathfrak
b}_\pm$ and with the condition $x_+|_{\mathfrak
h}=-x_-|_{\mathfrak h}$, where ${\mathfrak h}$ is the
corresponding Cartan subalgebra.

\begin{rem} We can choose
$G_+=\{e\}$ and $G_-=G$. We call it trivial factorizability.
\end{rem}

Let $G$ be a factorizable group. Define a binary
operation
\[
g\star h=g_+h_+(g_-h_-)^{-1}
\]
which obviously defines a group structure on $G$ with the
same identity element as for the original group structure. The
inverse to $g$ in this group is $i(g)=g_+^{-1}g_-$. For
factorizable Lie groups this structure is defined on a dense open
subset of $G$ (or on $G_+\times G_-$).

\subsection{Combinatorial description of the moduli space for factorizable
groups}\label{sec-fact}

Assume $G$ is a factorizable group and $D_t$ is a regular diagram of
a tangle.
\begin{defin}\label{fact-rel}
The map $E(D_t)\to G$ which associates to edge $e$
element $x_e\in G$ is called a $G$-coloring of
diagram $D_t$ if it satisfies the relations:
\begin{itemize}
\item double point:
$x_{b_v} = (x_{a_v})^{-1}_\pm x_{c_v} (x_{a_v})_\pm$,\
$x_{d_v} = (x_{c_v})_\mp^{-1} x_{a_v} (x_{c_v})_\mp$, for $(\pm)$ crossing,

\item critical point: $x_{e_v}=x_{r_v}$.
\end{itemize}
\end{defin}
\begin{thm} \label{color}
There is one to one correspondence between $G$-colorings of
a regular diagram of a tangle and
generic representations of the fundamental group of the
complement.
\end{thm}
\begin{proof}
We define a representation of
$\pi_1({\mathbb R}^2\times I\backslash t)$ which maps
the homotopy class of path $\Ga_e$ (see Fig.~\ref{mon})
to the element
\begin{equation} \label{bundling}
(x_{e_1})_-^{\ep_1}\dots (x_{e_n})_-^{\ep_n} \
x_e(x_{e_n})_-^{-\ep_n}\dots (x_{e_1})_-^{-\ep_1}
\end{equation}
where $e_1,\ldots,e_n$ is the sequence of edges, intersected by corresponding
to edge $e$
monodromy path (enumeration is given along the opposite direction
to its orientation), while
$\ep_i$ are the
signs defined by the projection of edges $e_i$ to vertical axis:
\[
\ep_i=\left\{ \begin{array}{rl} 1 & \mbox{if $e_i$ is oriented
upwards}\\ -1 & \mbox{if $e_i$ is oriented downwards}
\end{array}
\right .
\]
In other words, the holonomy along a monodromy path is computed as a
product of elements corresponding to "elementary intersections"
between edges of the diagram and this path, taken in the order
opposite to the orientation of the path.
%\begin{figure}
%\centerline{\psfig{figure=fig-holonom-rule.eps,width=1.2in}}
%\caption{\label{holonom-rule}}
%\end{figure}
Elements corresponding
to the "elementary intersections" are given on Fig.~\ref{crossing-rule}.
\begin{figure}
\centerline{\psfig{figure=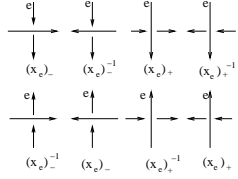,width=1.2in}}
\caption{\label{crossing-rule}There horizontal lines are parts
of the monodromy path while vertical lines are
edges of the diagram of a tangle.}
\end{figure}
To prove the theorem one has to check that holonomies
(\ref{bundling}) satisfy the relations defining the set $\Pi(D_t,G)$ and that
on an open dense subset of $\Pi(D_t, G)$ one can solve
equations (\ref{bundling}) recursively.
These are simple straightforward calculations.
\end{proof}

For factorizable $G$ we will make the standard choice of ramification
of a representation of the fundamental group of the complement.
Assume that the representation is generic and that $\{x_e\}$
is the $G$-coloring corresponding to it. Then the standard ramification
is obtained by applying rules from Fig. \ref{crossing-rule} to the
path ramification path connection points $x_0$ and $\tilde{x}_0$:
\[
g=(x_{e_1})_-^{-1}\dots (x_{e_n})_-^{-1}
\]
Here $e_1,\dots, e_n$ are edges of $D_t$ which intersect with the ramification
path.

\section{The category of tangles with flat connections in their
complements}

In this section we will define the category of tangles with
ramified representations of the fundamental group of the complement.
Since a flat connection defines a ramified representation of the
fundamental group, by the abuse of language we will call this category the
category of tangles with flat connections in the complement.

\subsection{The equivalence}

Restricted geometric tangle $t$ is called {\it
standard} if $\p_\pm t=\{(0,1,1),(0,2,1)\dots (0,n_\pm,1)\}$ for some
$n_\pm\in{\z}_{\ge0}$. We identify planes $\pl\x \{0\}$ and
$\pl\x\{1\}$ with $\pl$ through the projection $p$ defined in ~\eqref{eq:proj}.

Consider the set of pairs $(t,f,g)$, where $t$ is a standard tangle and
$f:\pi_1( \pl\times I\backslash t,x_0)\to G$ is a group homomorphism
and $g\in G$ is its ramification.
Define the following equivalence relation on this set.
Pair $(t,f,g)$ is equivalent to $(t',f',g)$
if $t$ is isotopic to $t'$ in ${\pl}\x I$ by an isotopy preserving
the boundary and points $x_0, \ \tilde{x}_0$ (so, in particular, $\p_{\pm}t =\p_{\pm}t'$),
which brings flat
connection $f$ to $f'$.
Equivalence class $[t,f,g]$ will be called a (standard) tangle with a
ramified representation of the fundamental group in the complement.

\subsection{The composition}\label{composition}
From now on the point $x_0$ will have coordinates
$(0,y_0,1/2)$ with $y_0<<-1$ and $\tilde{x}_0$ will
have coordinates $(0,\tilde{y}_0,1/2)$ with $\tilde{y}_0>>1$.

Let $t_1$ and $t_2$ be two  standard geometrical tangles with
$\p_+t_1=\p_-t_2$. Define
their composition $t_1\circ t_2$ as the tangle with the following
defining properties:
\begin{equation}\label{iso-1}
t_1\circ t_2\cap \pl\x [0,1/2]\simeq t_1,
\end{equation}
\begin{equation}\label{iso-2}
t_1\circ t_2\cap \pl\x [1/2,1]\simeq t_2.
\end{equation}
Here the identifications are given by linear maps,
$(x,y,z)\to (x,y,2z)$ for the first one, and
 $(x,y,z)\to(x,y,2z-1)$ for the second.

Let $[t_1,f_1,g]$ and $[t_2,f_2,g]$ be two tangles with flat
connection in their complements such that $\p_+ t_1=\p_- t_2$ and
$f_1\big\vert_{\pi_1({\pl}\backslash \p_+t_1,p(x_0))}
=f_2\big\vert_{\pi_1({\pl}\backslash
\p_-t_2,p(x_0))}$ and $f_1$ and $f_2$ have the same ramifications.
Define the mapping $f_{1\circ2}$
\[
 f_{1\circ2}\big\vert_{{\pl}\x [0,1/2]} = \tilde{f_1}, \quad
f_{1\circ2}\big\vert_{{\pl}\x [1/2,1]} = \tilde{f_2},
\]
where we identified elements of fundamental groups using
(\ref{iso-1}) and (\ref{iso-2}). The map $f_{1\circ2}$
is a representation of the fundamental group of the
complement of $t_1\circ t_2$.
We define the composition of two tangles with flat connections
$[t_1,f_1,g]$ and $[t_2,f_2,g]$ as the equivalence class
$[t_1\circ t_2,f_{1\circ2},g]$ and denote it as $[t_1,f_1,g]\circ[t_2,f_2,g]$.
It is clear that the composition of two tangles with flat connections
defined this way does not depend on the choice of
representatives in $[t_1,f_1,g]$, $[t_2,f_2,g]$, and it is associative.

\subsection{Tensor product}

Let $[t_1]$ and $[t_2]$ be isotopy classes of two standard tangles.
Let $|X|$ be the
cardinality of a finite set $X$. Consider the map $\var: \pl\x
I\to \pl\x I$ acting at the boundary as
\[
(x,y,0)\mapsto(x,  \ y+|\p_-t_1|,  \ 0)
\]
\[
(x,y,1)\mapsto(x,  \ y+|\p_+t_1|,  \ 1)
\]
Let $x_0(i)$ and $\tilde{x}_0(i)$ be the reference points for $t_i$
($i=1,2$).
Assume that the map $\var$ brings a
representative $t_2$ of $[t_2]$ to $t'_2$ and that it brings $x_0(2)$ to
$\tilde{x}_0(1)$. Choose representatives
of $t_1$ and $t_2$ in such a way that one dimensional submanifolds $t_1$
and $t'_2$ do not tangle with each other ( that their diagrams
do not intersect). The isotopy class of
$t_1\cup t'_2$  is called the tensor product $[t_1]\otimes [t_2]$
of $[t_1]$ and $[t_2]$ \cite{T-1}. It is clear that $[t_1]\otimes [t_2]$
does not depend on the choice of representatives.

Let $[t_1, f_1]$ and $[t_2,f_2]$ be two tangles with given representations
of the fundamental group of the complement. Let $t_1$ and $t'_2$ be geometrical tangles as above.
It is clear that $\pi_1(\pl\x I\backslash t_1\cup t'_2,x_0(1))$
is isomorphic to the group freely generated by $\pi_1(\pl\x I\backslash
t_1,x_0(1))$ and $\pi_1(\pl\x I\backslash t'_2, x_0(2))$.
Definition of the tensor product will depend on how we embed
groups $\pi_1(\pl\x I\backslash t_i,x_0(i))$ into
$\pi_1(\pl\x I\backslash t_1\cup t'_2, x_0(1))$.

\begin{figure}
\centerline{\psfig{figure=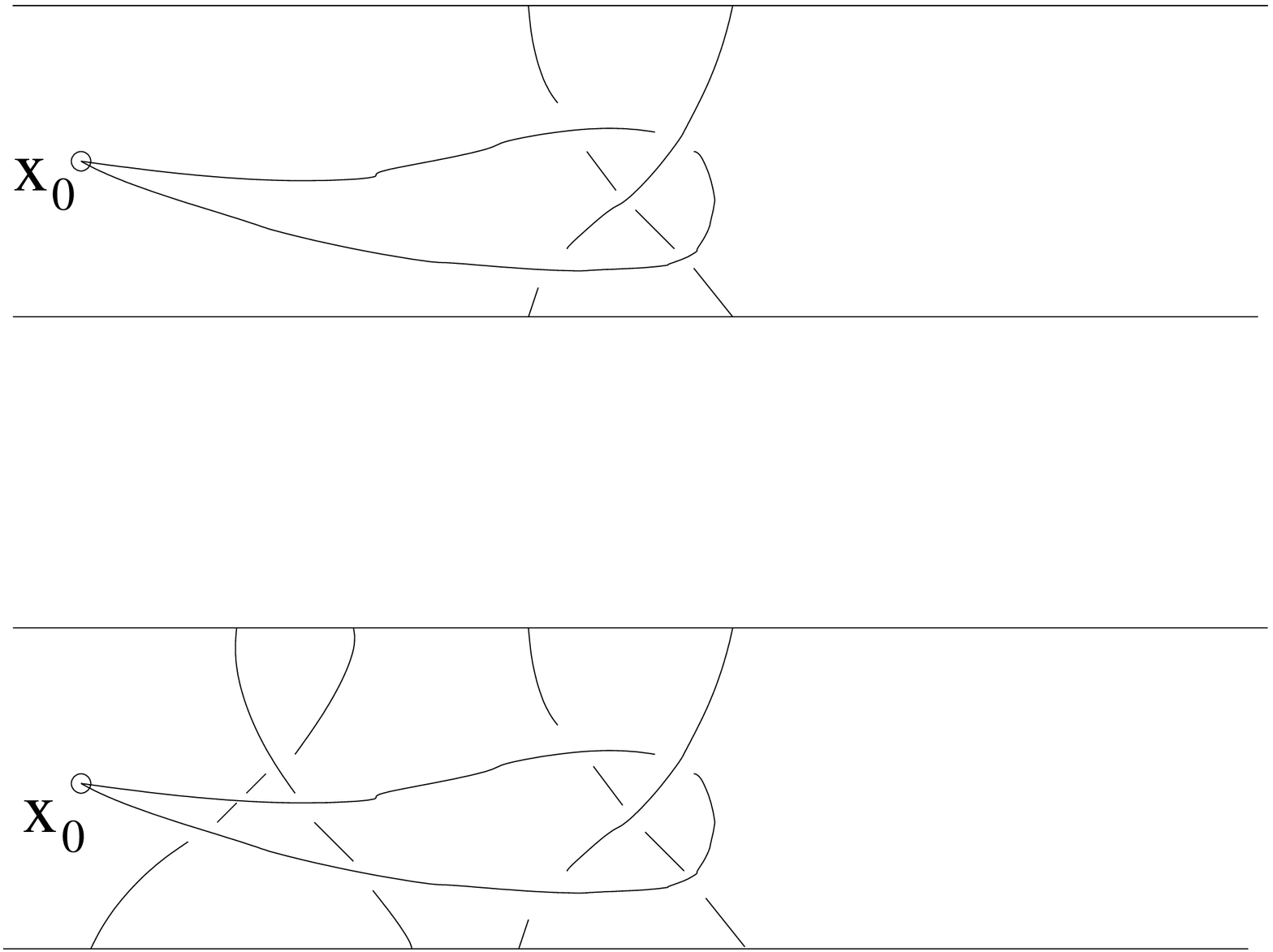,width=2.0in}}
\caption{\label{tp-of-fg}}
\end{figure}

Fix an isomorphism by identifying elements of
$\pi_1(\pl\x I\backslash t_1\cup t'_2,x_0(1))$ with corresponding elements of
$\pi_1(\pl\x I\backslash t_1,x_0(1))$ and
$\pi_1(\pl\x I\backslash t'_2, x_0(2))$
in accordance with Fig.~\ref{tp-of-fg}. This identification means that
we identify elements of $\pi_1(\pl\x I\backslash t'_2, x_0(2))$ with
elements of $\pi_1(\pl\x I\backslash t_1\cup t'_2,x_0(1))$ by gluing
them with the ramification path for $t_1$.
For a pair of group ramified
homomorphisms $f_i:\pi_1(\pl\x I\backslash
t_i,x_0)\to G$, $i=1,2$, with ramifications $g_1$ and $g_2$ respectively
define the group homomorphism
\[
f_{1\otimes_g2}\colon
\pi_1(\pl\x I\backslash
t_1\cup t'_2, x_0) \to G
\]
by assigning
\[
f_{1\otimes2}(\gamma)=\left\{\begin{array}{cc}
f_1(\gamma)&\mbox{if  }\gamma\in \pi_1(\pl\x I\backslash t_1,x_0)\\
g_1f_2(\gamma)g_1^{-1}&
\mbox{if  }\gamma\in \pi_1(\pl\x I\backslash t'_2,x'_0)
\end{array}\right.
\]
The ramification of $f_{1\otimes 2}$ we will define as $g_1g_2\in G$.
The equivalence class of $(t_1\cup t'_2, f_{1\otimes 2})$
will be called the
tensor product of $[t_1, f_1,g_1]$  and $[t_2,f_2,g_2]$ and will
be denoted $ [t_1, f_1,g_1]\otimes [t_2,f_2,g_2]$. This
operation is associative.

\subsection{Framing}A framing of a geometric tangle $t$ is a
continuous nonsingular section of the normal bundle to $t$.
A framing of a tangle is the homotopy class of a geometric framing.
We assume that the
framing is parallel to the $x$-axis at $\partial_\pm t$, has positive projection to this axis
and unit length.

\subsection{The category of framed G-tangles}
Here $G$-framed tangle is a tangle with a ramified representation of the
fundamental group of the complement.

Define the category of framed $G$-tangles $\widetilde{\mathcal FT}(G)$ as follows:
\begin{itemize}
\item {\it Objects} are  pairs $(\emptyset; f)$ where $f:\pi_1(\pl)\to e\in G$
with the trivial ramification ($g=1$)  and
$((\epsilon_1,\dots ,\epsilon_n);f,g)$ where
$\epsilon_i=\pm 1$, $f$ is a representation of the fundamental
group of ${\pl}\backslash\{(0,1),\dots ,(0,n)\}$ with the ramification $g$.

\item {\it Morphisms} between $((\epsilon_1,\dots ,\epsilon_n);f_+;g_+)$
and $((\sigma_1,\dots,\sigma_k);f_-;g_-)$ are nontrivial only if
$g_+=g_-$. They are oriented framed tangles with flat connections
(with ramified representation of the fundamental group of the complement)
$[t,f,g]$ where $t$ is a tangle with $|\p_-t|=k$, $|\p_+t|=n$ and
$f$ is a flat connection in the complement of $t$ such that
$f|_{\pi_1(\pl\backslash \p_\pm t, x_0)}=f_\pm$. The ramification $g$ is defined
by the ramifications of the objects $g=g_+=g_-$. The orientation
of boundary components should agree with signs
$\epsilon_i,\sigma_j$ as is shown on Fig.~\ref{directions}.

\begin{figure}
\centerline{\psfig{figure=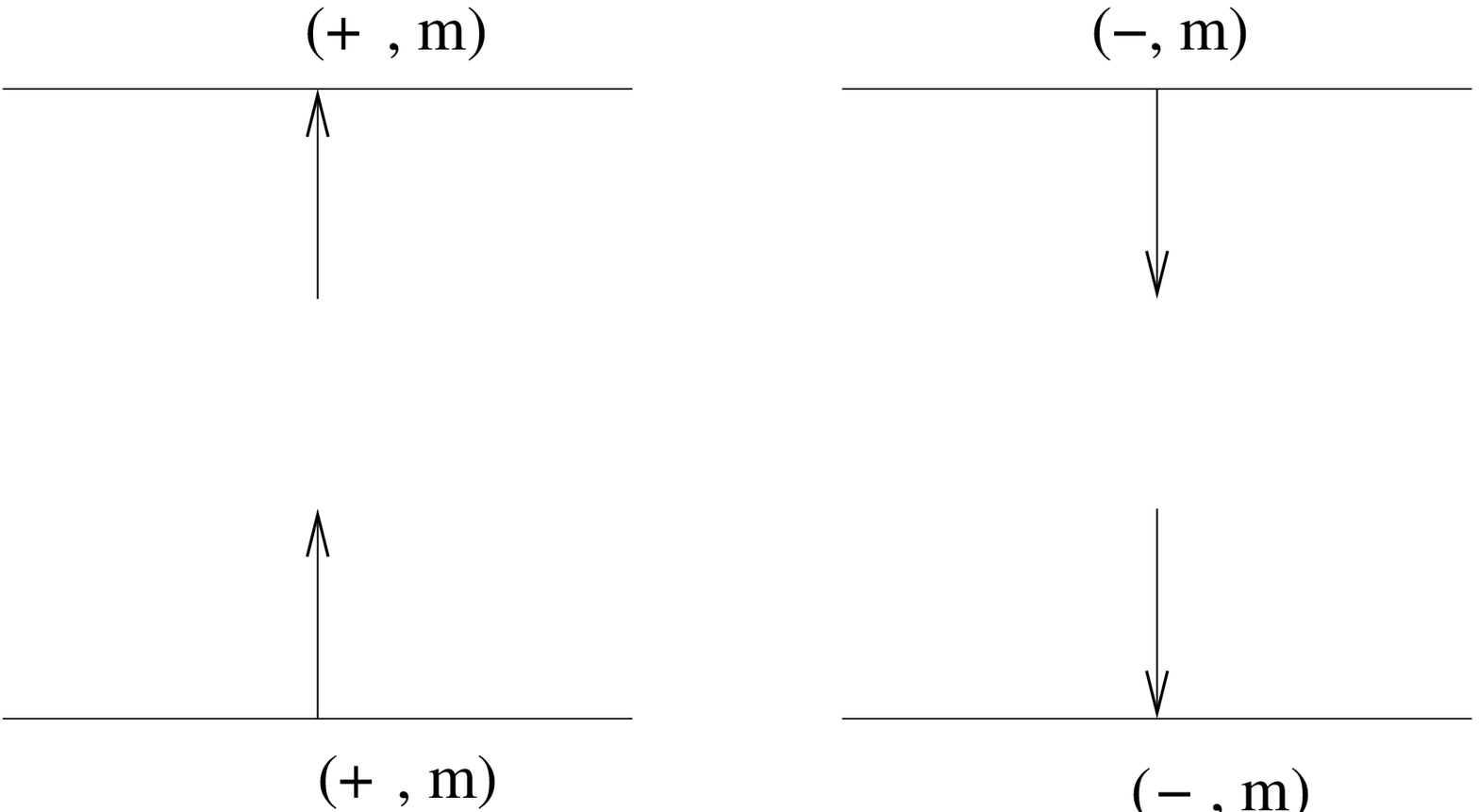,width=2.0in}}
\caption{\label{directions}}
\end{figure}

\item The composition is defined above in Subsection~\ref{composition}.

\item The identity element in End$((\epsilon_1,\dots,\epsilon_n),f,g)$
consists of the identity tangle (see Fig.~\ref{identity}) with the
representation of the fundamental group and the ramification uniquely defined by its
values on the boundary.
\end{itemize}

\begin{figure}
\centerline{\psfig{figure=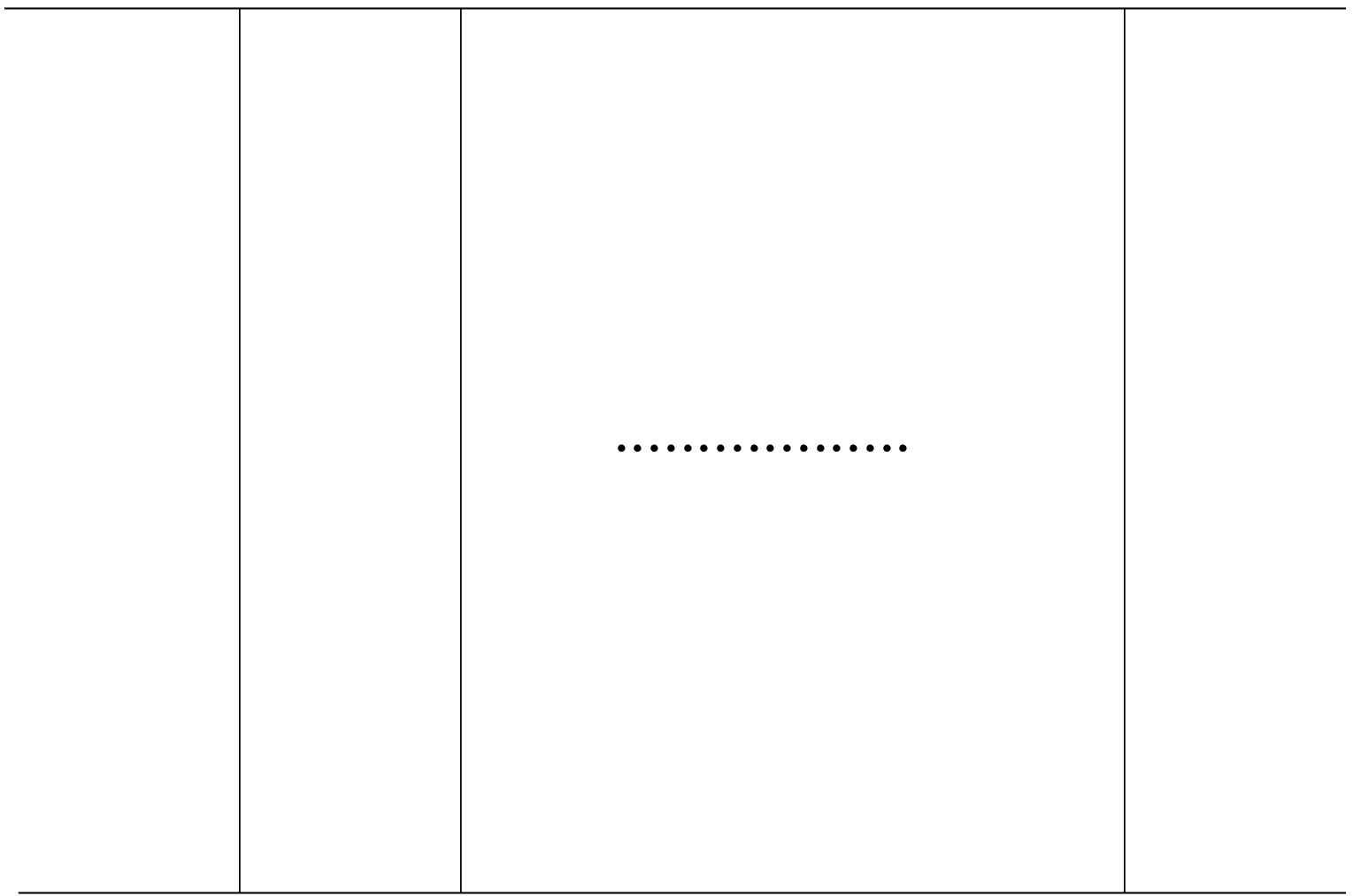,width=2.0in}}
\caption{\label{identity}}
\end{figure}

The category $\widetilde{\mathcal FT}(G)$ is braided monoidal with the
tensor product described above.

%The dual object to $((\epsilon_1,\dots ,\epsilon_n);f)$ is
%\[
%((\epsilon_1,\dots,\epsilon_n);f)^*\equiv((-\epsilon_n,\dots,-\epsilon_1)
%;\bar{f})
%\]
%where $\bar{f}$ is obtained from $f$
%by the reflection with respect to the plane
%$(x,y,z)=(\frac{n+1}{2},y,z)$.
%Diagrams of tangles representing
%evaluation and injection morphisms are given in Figures~\ref{evaluation} and
%\ref{injection}.

%\begin{figure}
%\centerline{\psfig{figure=fig-12.eps,width=2.0in}}
%\caption{\label{injection}}
%\end{figure}

%\begin{figure}
%\centerline{\psfig{figure=fig-13.eps,width=2.0in}}
%\caption{\label{evaluation}}
%\end{figure}

\begin{figure}
\centerline{\psfig{figure=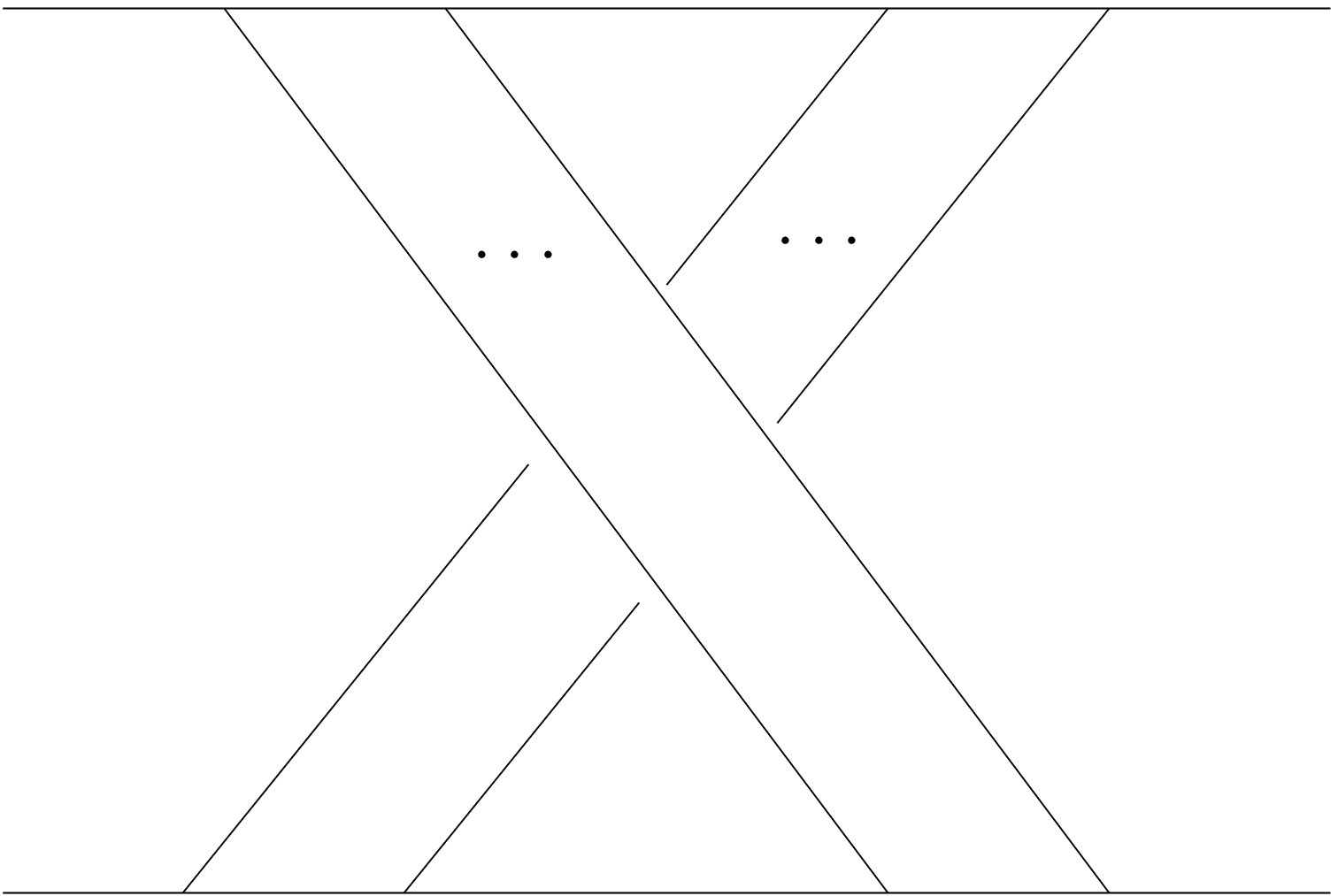,width=2.0in}}
\caption{\label{braiding}}
\end{figure}

The braiding is given by the tangle corresponding to the diagram
from Fig. \ref{braiding}. The corresponding ramified representation of the
fundamental group of its complement is completely determined by the
boundary values which are fixed by the choice of objects.

\subsection{The combinatorial description of $\widetilde{\mathcal FT}(G)$}

The category $\widetilde{\mathcal FT}(G)$ can be described
in terms of equivalence classes of diagrams with elements of $G$
assigned to edges of the diagram.

\begin{figure}
\centerline{\psfig{figure=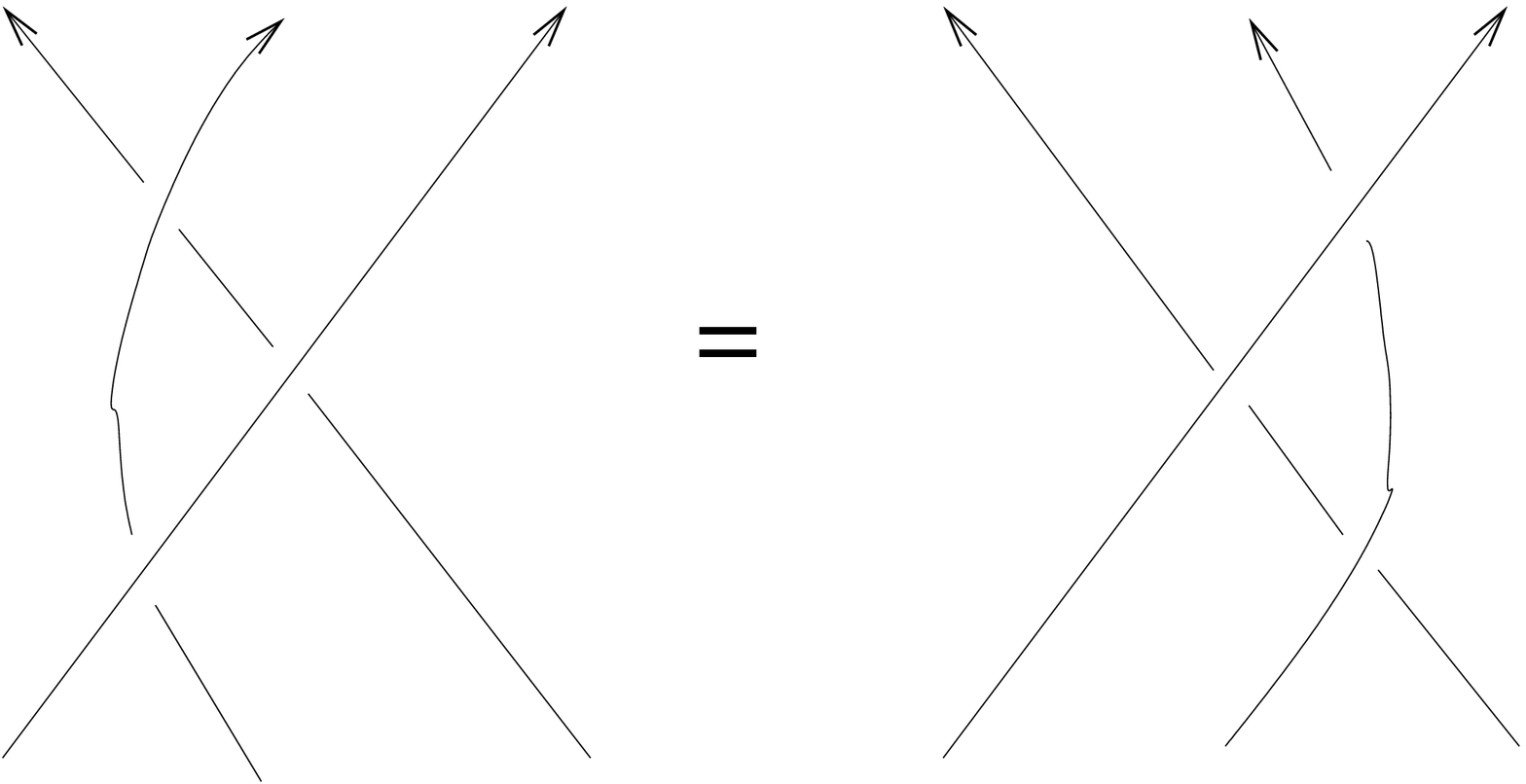,width=2.0in}}
\caption{\label{Red-1}}
\end{figure}

\begin{figure}
\centerline{\psfig{figure=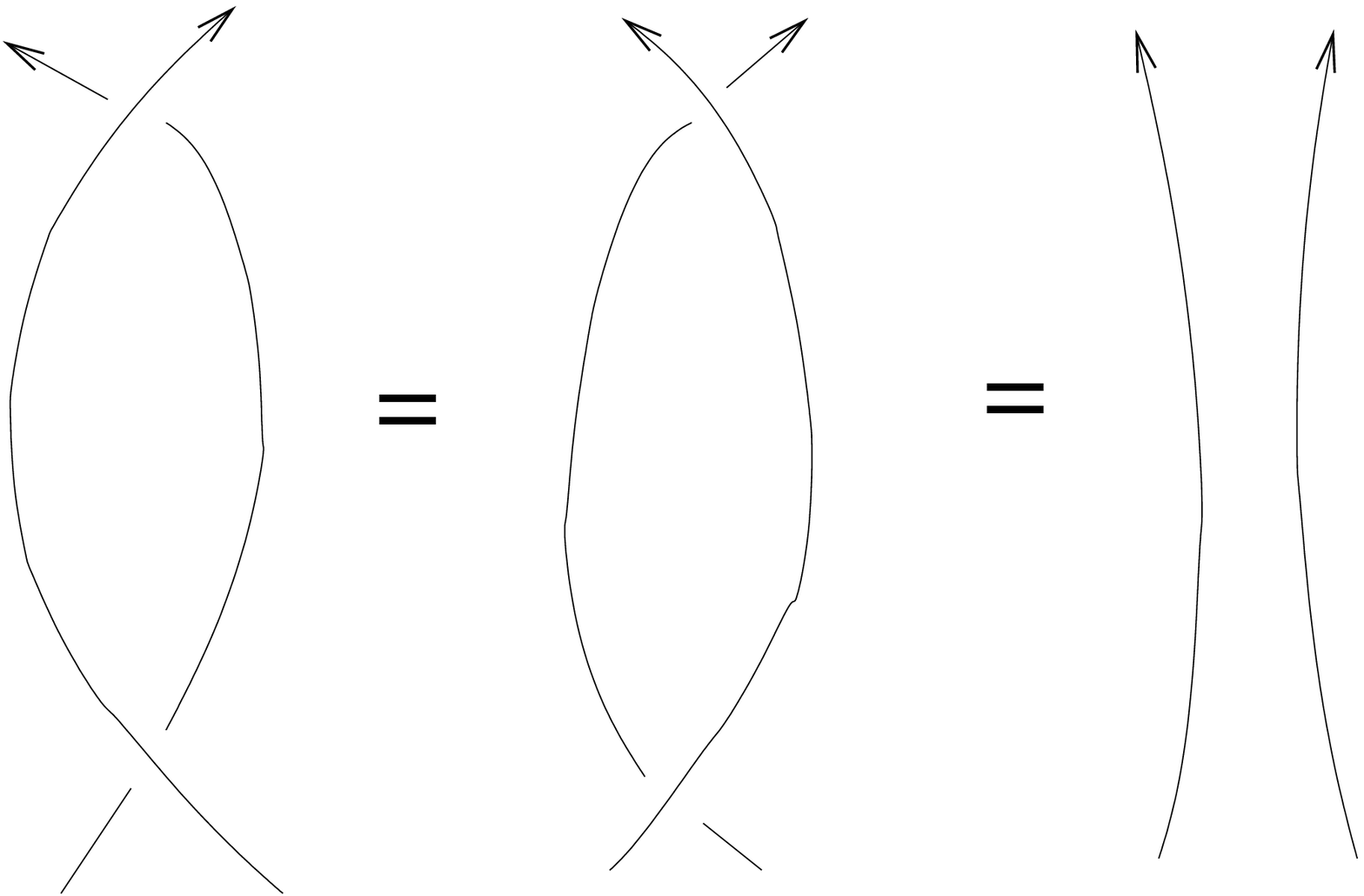,width=2.0in}}
\caption{\label{Red-2}}
\end{figure}

\begin{figure}
\centerline{\psfig{figure=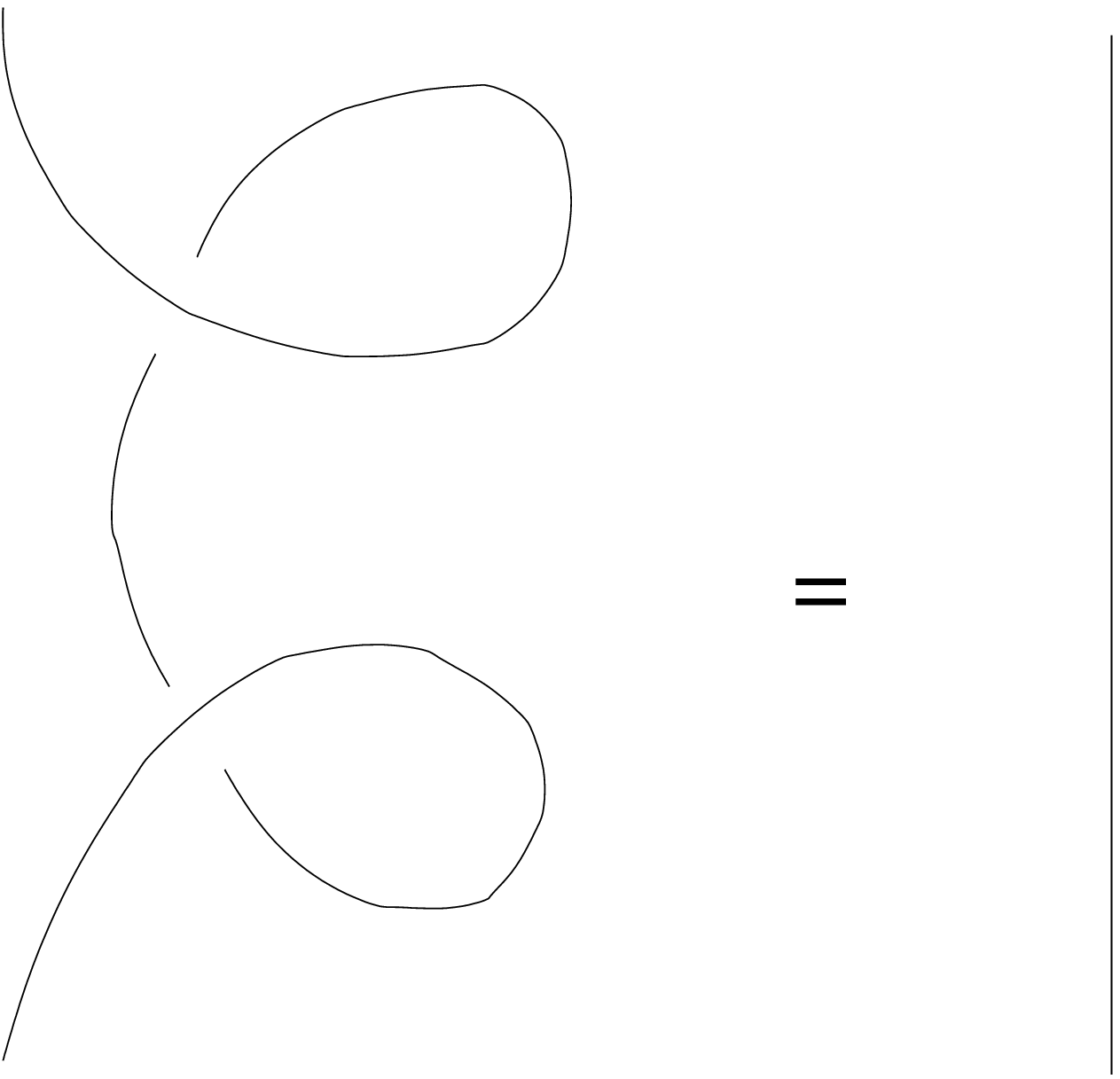,width=2.0in}}
\caption{\label{Red-3}}
\end{figure}

Representations of the fundamental group of the complement
to a tangle are in bijections with maps $E(D_t)\to
G$ which satisfy Wirtinger relations. Therefore a tangle with a flat
connection in the complement can be regarded as an equivalence
class of diagrams equipped with such data modulo Redemeister moves. If the tangle is
framed, one has to consider only framed Reidemeister moves (see
Figures~\ref{Red-1}, \ref{Red-2} \ref{Red-3}).

\begin{figure}
\centerline{\psfig{figure=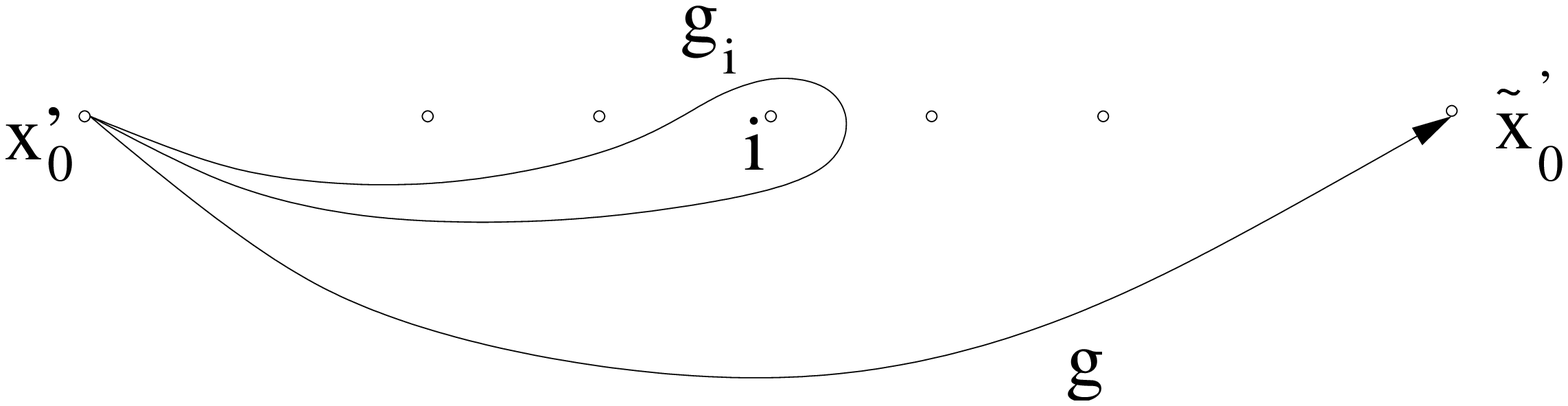,width=2.0in}}
\caption{\label{fund-cycles}}
\end{figure}

A representation of the fundamental group of the complement to $n$ distinct points in
$\pl$ is fixed uniquely by the monodromies along a system of
fundamental cycles. Choose them as on Fig.~\ref{fund-cycles} then
an object of the category $\widetilde{\mathcal FT}(G)$ can be
written as a collection $(\epsilon_1,\dots,\epsilon_n; \
g_1,\dots, g_n; \ g)$ where $g_i\in G$ is a monodromy along the cycle
shown on Fig.~\ref{fund-cycles} around the point $(0,i,0)$ and $g$ is
the ramification.

Thus, we can describe the category $\widetilde{\mathcal FT}(G)$
as follows:
\begin{itemize}
\item {\it Objects} are collections $(\ep_1,\dots,\ep_n;
g_1,\dots,g_n;g)$,
where $\ep_i=\pm$ and $g_i,g\in G$.

\item {\it Morphisms} between the object $(\ep_1,\dots,\ep_n; g_1,\dots,g_n;g)$
and the object $(\sigma_1,\dots,\sigma_m;h_1,\dots, h_m;g)\}$ are diagrams
of oriented framed tangles with elements of $G$ assigned to edges modulo Redemeister
moves and Wirtinger relations. The orientation and group elements
assigned to boundary components should agree the orientation and group
elements of corresponding boundary components.

\item {\it Composition} is the usual composition of diagrams of tangles
with group elements assigned to edges.

\item {\it Tensor product} for objects is $(\{\ep_i\};\{g_i\};g)\otimes
(\{\sigma_i\};\{h_i\};h)=(\{\ep_i\}\{\sigma_i\};\{g_i\}\{h_i\};gh)$.
The tensor product of two morphisms $(D_1,\{g_e\},g)$ and $(D_2,\{h_e\},h$
is the diagram which $D_1\sqcup D_2$ is
the tensor product of two diagrams (juxtaposition) with the elements $g_e$
assigned to the edges of $D_1$ and elements $g^{-1}h_eg$ assigned to edges
of $D_2$. The ramification of the tensor product is $gh$.

\end{itemize}

\subsection{Combinatorial description of $\widetilde{\mathcal
FT}(G)$ for factorizable $G$} Let $G$ be a factorizable group.
According to Theorem~\ref{color} a $G$-coloring of diagram
$D_t$ is in one to one correspondence with a flat $G$-connection
in the complement to
tangle $t$ by assigning monodromies $\ga_e$ to edges
of the diagram:
\begin{equation}\label{themap}
\{x_e\}_{e\in D_t} \mapsto  \{\ga_e=(x_{e_1})_-^{\ep_1}\dots
(x_{e_n})_-^{\ep_n} \ x_e(x_{e_n})_-^{-\ep_n}\dots
(x_{e_1})_-^{-\ep_1}\}_{e\in D_t}
\end{equation}
More precisely, it is a
bijection for factorizable groups and a bijection on an open dense
subset of $G^{\times E(D_t)}$  for almost factorizable Lie groups.
This suggests the description of the the category $\widetilde{\mathcal
FT}(G)$ in term of $G$-colored diagrams.

Define the category ${\mathcal FT}(G)$ of $G$-colored
framed diagrams as follows:

\begin{itemize}
\item {\it Objects} are sequences $\{(\epsilon_1,x_1),\dots,
(\epsilon_n,x_n)\}$ where $\e_i=\pm 1$ and
generic $x_i\in G$.
\item {\it Morphisms} are $G$-colored framed diagrams modulo framed
Reidemeister moves (see Figures~\ref{Red-1}, \ref{Red-2}, \ref{Red-3})
\end{itemize}

This category is monoidal, braided. The composition is
obvious. The tensor product is essentially the same as the tensor
product of diagrams of framed tangles (see for example \cite{T-1}).
The decoration of the tensor product is defined by the decoration
of factors.

%The dual object to
%$\{(\epsilon_1,x_1),\dots,(\epsilon_n,x_n)\}$ is
%$\{(-\epsilon_n,x_n),\dots,(-\epsilon_1,x_1)\}$.

Let $\widetilde{\mathcal FT}$ be the category of tangles with flat
connections in the complement with the standard ramification
(see section \ref{sec-fact}).

\begin{prop}Consider the map  $E: {\mathcal FT}(G)\to \widetilde{\mathcal FT}(G)$
acting on objects as
\[
\{ (\epsilon_1,x_1),\dots, (\epsilon_n,x_n)\}\to
\{(\epsilon_1,g_1),\dots ,(\epsilon_n,g_n)\}
\]
where
\[
g_i=(x_1)_-^{\ep_1}\dots
(x_{i-1})_-^{\ep_{i-1}}x_i
(x_{i-1})_-^{-\ep_{i-1}}\dots(x_{1})_-^{-\ep_1}
\]
and on morphisms as $(D_t$ with  a $G$-coloring$)\to (t$ with
ramified representations of the $\pi_1$ of the complement to $t$ defined by
the $G$-coloring of $D_t)$. It is a monoidal, braided functor
${\mathcal FT}(G) \to \widetilde{{\mathcal FT}}(G)$. This
functor is also an equivalence of categories.
\end{prop}

The proof follows from the bijection between tangles with
given representation of the fundamental group of the complement and
the equivalence classes of $G$-colored diagrams.

Because the categories $\widetilde{\mathcal FT}(G)$and ${\mathcal FT}(G)$
are naturally equivalent we will abuse notations and will denote
both of them by ${\mathcal FT}(G)$.

When $G$ is an almost factorizable Lie group the category of
$G$-colored diagrams has the tensor product, but the dual object
exists only on an open dense subset of objects. The braiding in
this case also exists only on an open subset of objects of this
category.

\subsection{Elementary diagrams}

The following fact is a key for construction of invariants of
tangles via braided monoidal categories.

\begin{figure}
\centerline{\psfig{figure=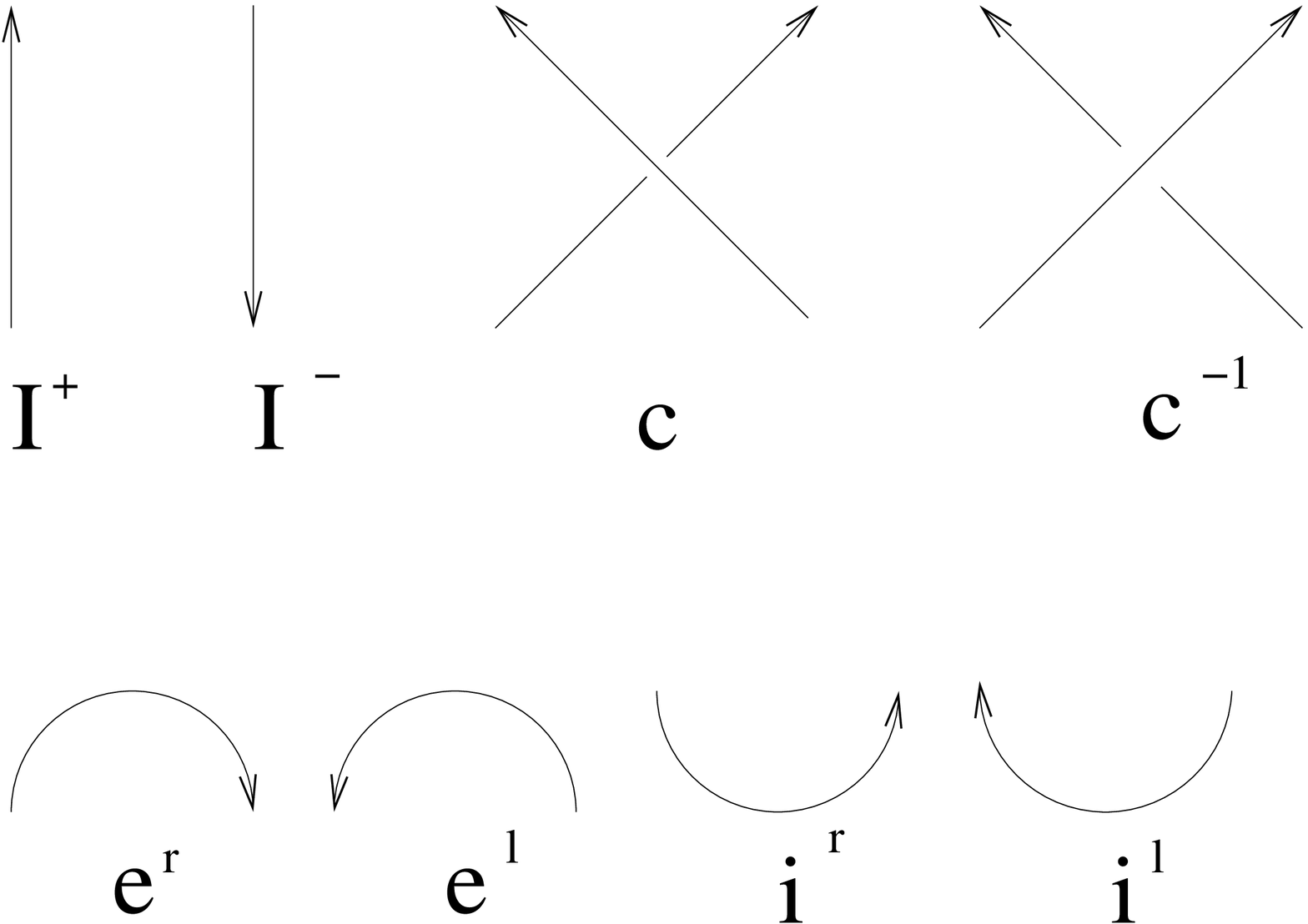,width=2.0in}}
\caption{\label{elem-diagr}}
\end{figure}

\begin{prop}
Framed $G$-colored tangles are compositions of tensor products of
elementary diagrams. Elementary diagrams are given on Fig.
\ref{elem-diagr} .
\end{prop}

The proof of this proposition and the definition of elementary
diagrams are the same as for the category of framed tangles (see
\cite{T-1}).

\section{Holonomy Yang-Baxter equation and invariants of tangles with
flat connections}

Let $G$ be a factorizable group (or a factorizable Lie group).

\subsection{Holonomy Yang-Baxter equation}

Assume we have the following data:

\begin{itemize}
\item Collection of finite-dimensional
vector spaces ${\mathfrak X}$
\item For each pair $(X,Y)$ of vector spaces from this collection
we have a map
\begin{equation}\label{R}
R^{X, Y}: G\times G \to End(X\otimes Y)
\end{equation}
(or, more generally, a section of a  vector bundle over
$G\times G$ with the fiber $X\otimes Y$).
\end{itemize}

Let $x_L, x_R: G\times G\to G$ be mappings acting as
\begin{equation}\label{xLR}
x_L(x,y)=x_-yx_-^{-1}, \ x_R(x,y)=x_L(x,y)_+^{-1}x x_L(x,y)_+
\end{equation}

\begin{defin} The data described above define a system of holonomy
$R$-matrices if the equation
\begin{eqnarray}\label{hybe}
R^{X,Y}_{12}(x_R(x,x_L(y,z)),x_R(y,z)) R^{X,Z}_{13}(x,x_L(y,z))
R^{Y,Z}_{23}(y,z)=
\\ R^{Y,Z}_{23}(x_L(x,y),x_L(x_R(x,y),z)) R^{X,Z}_{13}(x_R(x,y),z) R^{X,Y}_{12}(x,y)
\end{eqnarray}
holds for any generic triple $x,y,z\in G$ and any $X,Y,Z\in
\mathfrak X$. Here the mappings $x_{L,R}$ are as above, all
factors act in $X\otimes Y\otimes Z$ and lower indices in the
$R$-matrices indicate in which factors of the tensor product they
act nontrivially.
\end{defin}

We will call the equation (\ref{hybe}) the {\it holonomy
Yang-Baxter equation}.

For a liner operator $f:V\to V$ we denote by $f^t$ the dual linear
operator $f^t:V^*\to V^*$.

\begin{defin} Let $V$ and $W$ be two finite-dimensional
vector spaces and $W^*$ be the vector space dual to $W$. We will
say that the liner map $R: V\otimes W\to V\otimes W$ is a {\emph
cross-nondegenerate } if the corresponding linear map
$R^{t_2}:V\otimes W^*\to V\otimes W^*$ is nondegenerate.
\end{defin}

For cross-nondegenerate holonomy $R$-matrix acting in $X\otimes X$
define the map ( more generally a section of the trivial vector
bundle) $d^X: G\to End(X)$.

\[
d^X(a)=(tr\otimes id)(P((R^{t_1}(a,i(a)^{-1})^{-1})^{t_1})
\]
Here and below we will use the notation $t_i$ for taking dual to
an operator acting in a tensor product with respect to the $i$-th
factor.

\begin{lemma} For generic $a\in G$ the linear operator $d^X(a)$
is invertible and
\[
d^X(a)^{-1}=(id\otimes
tr)((((R(i(a)^{-1},a)^{-1})^{t_2})^{-1})^{t_2}P)
\]
\end{lemma}

\begin{prop} For $X,Y\in \mathfrak X$ and generic $x,y\in G$ the
following identities hold
\begin{equation}\label{cr-1}
R_{12}^{X,Y}(x,y)=d_2(a)^{-1}(((R_{12}^{X,Y}(x,y)
^{-1})^{t_2})^{-1})^{t_2})d_2(y)
\end{equation}
\begin{equation}\label{cr-2}
R_{12}^{X,Y}(a,b)=d_1(y)^{-1}(((R_{12}^{X,Y}(a,b)
^{t_1})^{-1})^{t_1})^{-1})d_1(a)
\end{equation}
\begin{equation}\label{comm}
d_1^X(b)d_2^Y(a)R^{X,Y}(x,y)=R^{X,Y}(x,y)d_1^X(x)d_2^Y(y)
\end{equation}
where $a=x_L(x,y)$ and $b=x_R(x,y)$.
\end{prop}

For $X\in \mathfrak X$ and generic $x\in G$ define the linear
operator $w_X(x)\in End(X)$ as
\begin{equation}\label{w}
w_X(x)=tr_1 tr_2
(P_{12}((R_{12}^{-1}(x,b)^{t_2})^{-1})^{t_2}P_{13}R_{13}(b,x))
\end{equation}
where $b=i(x^{-1})$.

\begin{prop} For each $X,Y\in \mathfrak X$ we have
\[
(1\otimes w_Y(x_L(x,y))R^{X,Y}(x,y)=R^{X,Y}(x,y)(1\otimes w_Y(y))
\]
\[
(w_X(x_R(x,y))\otimes 1)R^{X,Y}(x,y)=R^{X,Y}(x,y)(w_X(x)\otimes 1)
\]
\end{prop}

The proof of this proposition and of the lemma are very similar.
Each of the equalities follows from a sequence of equalities
where each one corresponds to an elementary
moves of corresponding diagrams. They are completely parallel
to the proofs of similar statements from \cite{R-1}
( see Fig. 8, Fig. 13-17 from that paper and corresponding sequences of
equalities). The only difference is that $G$-coloring should be taken into
account.

\subsection{The category ${\mathcal FT}(G,{\mathfrak X})$}

Let $(D, \{x_e\}_{e\in E(D_t)})$ be a $G$-colored diagram of a
tangle. We will call it {\it decorated by  vector spaces
$\mathfrak X$} if a vector space from $\mathfrak X$ is assigned to
each connected component of the corresponding tangle.

Define the category ${\mathcal FC}(G,{\mathfrak X})$ of
$\mathfrak X$-decorated $G$-colored diagrams as follows.
\begin{itemize}
\item {\it Objects} of ${\mathcal FC}(G,{\mathfrak X})$ are sequences
$\{(\e_1,g_1,X_1), \dots (\e_n, g_n, X_n)\}$ where
$\e_i=\pm 1$, $g_i\in G'\subset G$ (recall that
$G'\subset G$ is the set of elements in $G$ which admit left and
right factorization) and $X_a\in {\mathfrak X}$.

\item {\it Morphisms} between two objects
$\{(\e_a,g_a, X_a)\}_{a=1}^n$ and
$\{(\sigma_a,h_a,Y_a)\}_{a=1}^k$ are $G$-colored diagrams
decorated by vector spaces from $\mathfrak X$ modulo framed
Redemeiseter moves. Decorations of boundary components are
determined by the initial and the target object. Elements $g_a$ and
$h_a$ are colorings of boundary edges of such diagrams.

\item Composition of morphisms is defined by gluing
decorated diagrams along the common boundary.

\end{itemize}

It is clear that $Ad_G$ of the $G$-coloring of a morphism in this category can be naturally
identified with a gauge class of a flat connections in the
complement to the corresponding tangle with connected components colored by $\mathfrak X$. The
category ${\mathcal FT}(G,{\mathfrak X})$ is braided monoidal. The
tensor product and the braiding are essentially the same as in the
category ${\mathcal FT}(G)$.
%The dual to $(\e,m,X,x)$ is $(-\e,X^*,i(x))$.

\subsection{The functor $F:
{\mathcal FT}(G,{\mathfrak X})\to
\underline{Vect}$}\label{invariant}

Here we will show how to use a system of holonomy $R$-matrices to
construct invariants of tangles with flat $G$-connections in the
complement. The construction is similar to the one of invariants
of tangles given in \cite{RT}. We will construct the monoidal
functor from the category $\tilde{{\mathcal F}C}(G,{\mathfrak X})$
to the category of all vector spaces.

\begin{thm}There exists a unique covariant functor $F:{\mathcal FT}
(G,{\mathfrak X})\to \underline{Vect}$ such that
\begin{itemize}
\item $F(\{(\e_1,g_1,X_1), \dots (\e_n,g_n, X_n)\})=
X_1^{\epsilon_1} \otimes \dots \otimes X_n^{\epsilon_n}$,
where $X^+=X, \ X^-=X^*$.

\item $F$ is a monoidal functor, i.e.
\[
F(\hat{D}_1\otimes \hat{D}_2)=F(\hat{D}_1)\otimes F(\hat{D}_2)
\]
where $\hat{D}_i$ are decorated diagrams.
\item Values of $F$ on elementary diagrams are:
\begin{enumerate}
\item If $I^\e: (\e,X,x)\to (\e,X,x)$ is the identity morphism
then
\[
F(I^\e)=id_{X^\epsilon}.
\]

\item For the morphism  $e^r: (\e,X,x)^*\otimes (\e,X,x)\to 1$ we have:
\[
F(e^r)=\left\{ \begin{array}{ll} e_X: X^*\otimes X\to \C,
l\otimes  v \mapsto l(d_X(x)v), & \mbox{if $\e=+1$} \\
X\otimes X^*\to \C, v\otimes l\mapsto l(v), & \mbox{if $\e=-1$}
\end{array} \right .
\]

\item For $i^r: 1\to(\e,X,x)\otimes (\e,X,x)^*$ we have
\[
F(i^r)=\left\{ \begin{array}{ll} i_{X}: \C\to X\otimes X^* ,
1\mapsto \sum_i d_X(x)^{-1}e_i\otimes e^i, & \mbox{ if $\e=+1$} \\
\C\to X^*\otimes X, 1\mapsto \sum_i e^i\otimes e_i, & \mbox{ if $\e=-1$}
\end{array}\right .
\]
Here $e_i$ is a linear basis in $X$ and $e^i$ is the
dual linear basis in $X^*$.
\item If $c: (+,X,x)\otimes (+,Y,y)\to (+,Y,y)\otimes (+,X,x)$ is the
braiding morphism,
\[
F(c)=P^{X,Y}R^{X,Y}(x,y): \ \ X\otimes Y \to Y\otimes X.
\]
\[
F(c^{-1})= (P^{X,Y}R^{X,Y)}(x,y))^{-1}: \ \ X\otimes Y \to
Y\otimes X.
\]
Here $P^{V,W}: V\otimes W\to W\otimes V$ is the linear map
$P^{V,W}(x\otimes y)= y\otimes x$.

\end{enumerate}

%\item The functor $F$ is rigid and braided.
\end{itemize}
\end{thm}

The proof of this theorem is completely parallel to the proof from of
the corresponding theorem for
tangles without flat connections (see \cite{RT}). One should
check the defining relations between elementary tangles, which is
a routine exercise.

The value of $F$ on a diagram is a linear map between the vector
spaces which are defined by the decoration of the boundary of the
diagram. The value of the functor $F$ on a diagram depends only on
its isotopy class and therefore $F$ is an invariant of framed tangles.
In the next section we will show that the functor $F$ is gauge invariant.

\section{Gauge invariance}

\subsection{The gauge group action and the dressing action}

\begin{figure}
\centerline{\psfig{figure=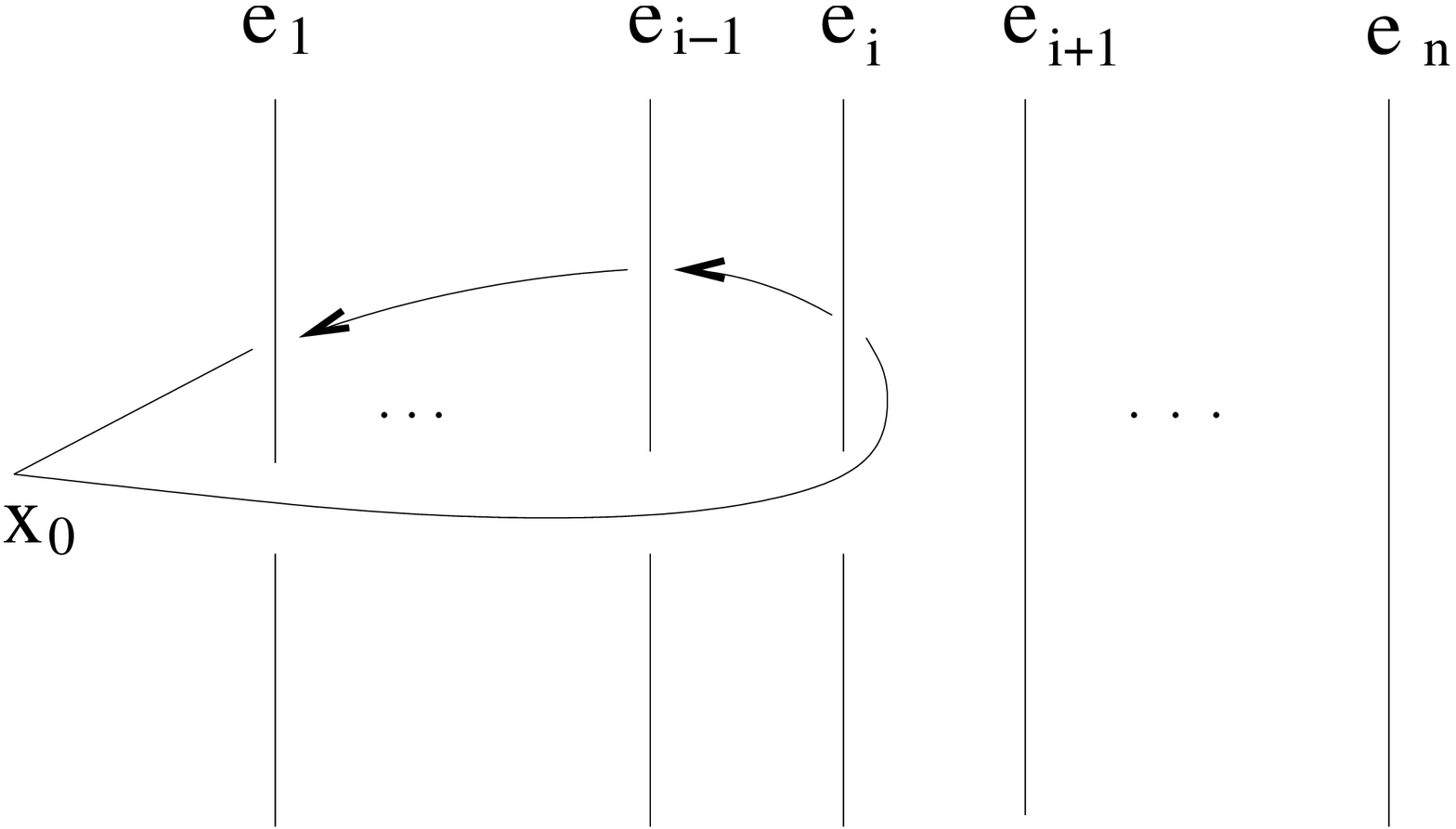,width=1.2in}}
\caption{\label{holonom}}
\end{figure}

Let $g_1, \dots g_n$ be a $G$-coloring of edges $e_1, \dots, e_n$
on Fig. \ref{holonom} and let $\gamma_i$ be the holonomy along the
path $\Gamma_i$. For generic colorings we have:
\begin{equation}\label{g-gamma}
\gamma_i=(g_1)_+^{\ep_1}\dots
(g_i)_+^{\ep_i}(g_i)_-^{-\ep_i}\dots(g_1)_-^{-\ep_1}
\end{equation}
or, equivalently,
\begin{equation}\label{gamma-g}
\ep_i(g_i)=(\gamma_{i-1})_+^{-1}\gamma_i(\gamma_{i-1})_-
\end{equation}

Here $\ep_i$ is the orientation of the edge $e_i$ ($+$ for going
up and $-$ for going down) and
\[
\ep(g)=\left\{ \begin{array}{cc} g & \mbox{for $\ep=+1$} \\ i(g) &
\mbox{for $\ep= -1$}
\end{array}
\right.
\]

Holonomies $\gamma_i$ and colorings $g_i$ can be regarded as an
element of $G^{\times n}$. Then formulae (\ref{g-gamma}),
(\ref{gamma-g}) describe two mutually inverse birational mappings
of $G^{\times n}$ to itself.

The gauge group acts on $\gamma_i$ by diagonal conjugation in
$G^{\times n}$:
\begin{equation}\label{holon-map}
x:\gamma_i\mapsto x\gamma_i x^{-1}
\end{equation}

The gauge action of $G$ on this space is the diagonal action by
conjugations.

For $x\in G$ consider the (birational)  mapping $\alpha_x:
G^{\times n}\to G^{\times n}$ defined as:
\begin{equation}\label{color-map}
\ep_i(\alpha_x(g_i))= x^{\gamma_{i-1}}\ep_i(g_i)(x^{\gamma_{i-1}})^{-1}
\end{equation}
where $\gamma_i$ are as above and
\begin{equation}\label{dress}
x^\gamma=(x\gamma x^{-1})_+^{-1}x \gamma_+= (x\gamma
x^{-1})_-^{-1}x \gamma_-
\end{equation}
is the {\it dressing} action of $\gamma$ on $x$.

Notice that the dressing action preserves the subgroups $G_\pm$
and
\begin{equation}\label{+dress}
(x_+)^\gamma=(\gamma_-^{-1}x_+^{-1})_+^{-1}
\end{equation}
\begin{equation}\label{-dress}
(x_-)^\gamma=(x_-\gamma_+)_-^{-1}
\end{equation}

\begin{prop} The mapping (\ref{gamma-g}) intertwines the gauge
action of $x\in G$ and the mapping $\alpha_x$.
\end{prop}

\begin{proof} Let $g^x=\{ (g^x)_1,\dots, (g^x)_n \}$ be the result of the
gauge action of $x\in G$ on $g$-coordinates. We have:
\begin{multline*}
\ep_i((g^x)_i)=(x\gamma_{i-1}x^{-1})_+^{-1}x\gamma_i x^{-1}
(x\gamma_{i-1}x^{-1})_-\\
=(x\gamma_{i-1}x^{-1})_+^{-1}x(\gamma_{i-1})_+
\ep_i(g_i)(\gamma_{i-1})_-^{-1}x^{-1}(x\gamma_{i-1}x^{-1})_-=
x^{\gamma_{i-1}}\ep_i(g_i)(x^{\gamma_{i-1}})^{-1}=\ep_i(\alpha_x(g)_i).
\end{multline*}
Here we used formulae (\ref{dress}) for the dressing action.
\end{proof}

\begin{cor} If $x$ belongs to one of the subgroups $G_\pm\subset
G$ we have
\begin{equation}\label{alpha+}
\ep_i(\alpha_{x_+}(g)_i)=x_+^{(g_{i-1})_-^{-\ep_{i-1}}
\dots(g_{1})_-^{-\ep_{1}}}
\ep_i(g_i)\left(x_+^{(g_{i-1})_-^{-\ep_{i-1}}
\dots(g_{1})_-^{-\ep_{1}}}\right)^{-1}
\end{equation}
\begin{equation}\label{alpha-}
\ep_i(\alpha_{x_-}(g)_i)=x_-^{(g_1)_+^{\ep_1}\dots(g_{i-1})_+^{\ep_{i-1}}}
\ep_i(g_i)\left(x_-^{(g_1)_+^{\ep_1}\dots(g_{i-1})_+^{\ep_{i-1}}}\right)^{-1}
\end{equation}
\end{cor}

\begin{figure}
\centerline{\psfig{figure=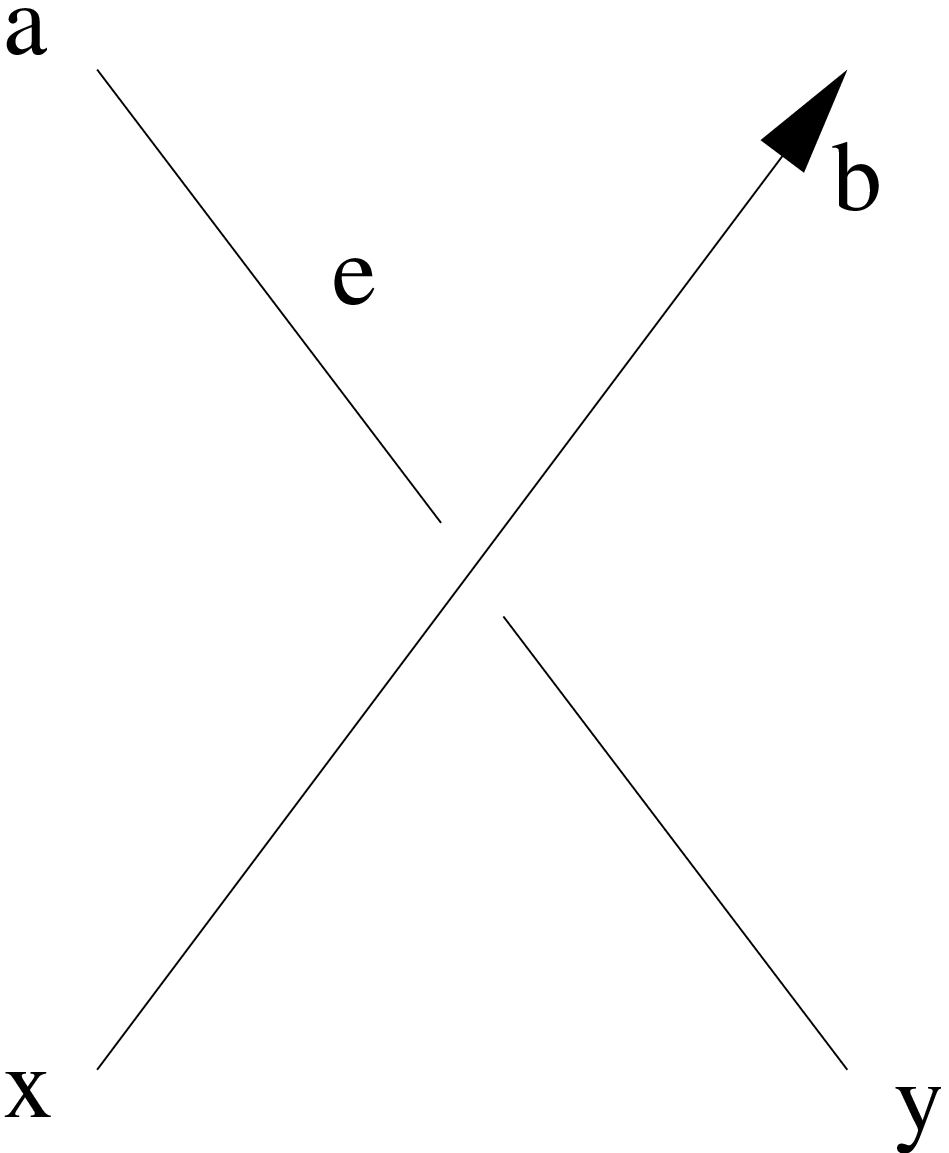,width=1.2in}}
\caption{\label{cross1}}
\end{figure}

\begin{lemma} \label{+cross} Let $x, y, a, b$ be $G$-colorings
of the diagram from Fig.~\ref{cross1}. Then,
\[
\ep_e(a)=x_+\ep_e(y)x_+^{-1}, \quad
b_+=\left(y_-^{-\epsilon_e}x_+^{-1}\right)_+^{-1}
\]
\end{lemma}

\begin{proof} According to the definition of $G$-coloring we have:
\[
b=a_-^{-\epsilon_e}xa_-^{\epsilon_e}, \quad
\epsilon_e(a)=x_+\epsilon_e(y)x_+^{-1}
\]
This gives the formula for $a$. In particular,
\[
a_-=\left(y_-^{-\epsilon_e}x_+^{-1}\right)_-^{\epsilon_e}
\]
so that
\[
b_+=\left(a_-^{-\epsilon_e}x_+\right)_+=
\left( \left(y_-^{-\epsilon_e}x_+^{-1}\right)_-^{-1}  x_+\right)_+
=\left( \left(y_-^{-\epsilon_e}x_+^{-1}\right)_+^{-1}
y_-^{-\epsilon_e}x_+^{-1}
x_+\right)_+=\left(y_-^{-\epsilon_e}x_+^{-1}\right)_+^{-1}
\]
\end{proof}

\begin{figure}
\centerline{\psfig{figure=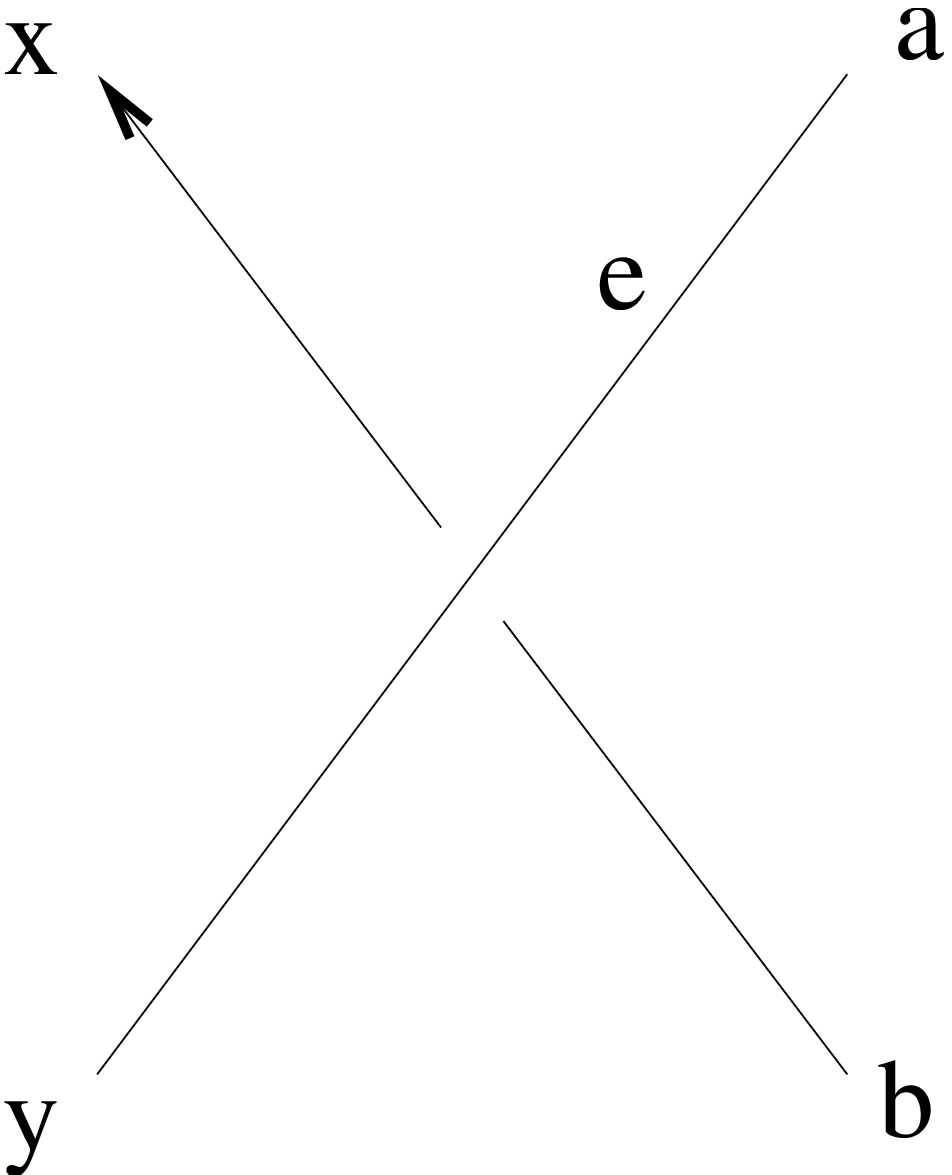,width=1.2in}}
\caption{\label{cross2}}
\end{figure}

\begin{lemma}\label{-cross} Let $x,y, a, b$ be colorings
the diagram from Fig.~\ref{cross2}. Then
\[
\epsilon_e(a)=x_-^{-1}\epsilon_e(y)x_-, \quad b_-=
\left(x_-^{-1}y_+^{\epsilon_e}\right)_-
\]
\end{lemma}

\begin{proof} From the definition of $G$-coloring we have
\[
\epsilon_e(a)=x_-^{-1}\epsilon_e(y)x_-,\quad
b=y_+^{-\epsilon_e}xy_+^{\epsilon_e}
\]
so the required formulas are straightforward.
\end{proof}

\begin{figure}
\centerline{\psfig{figure=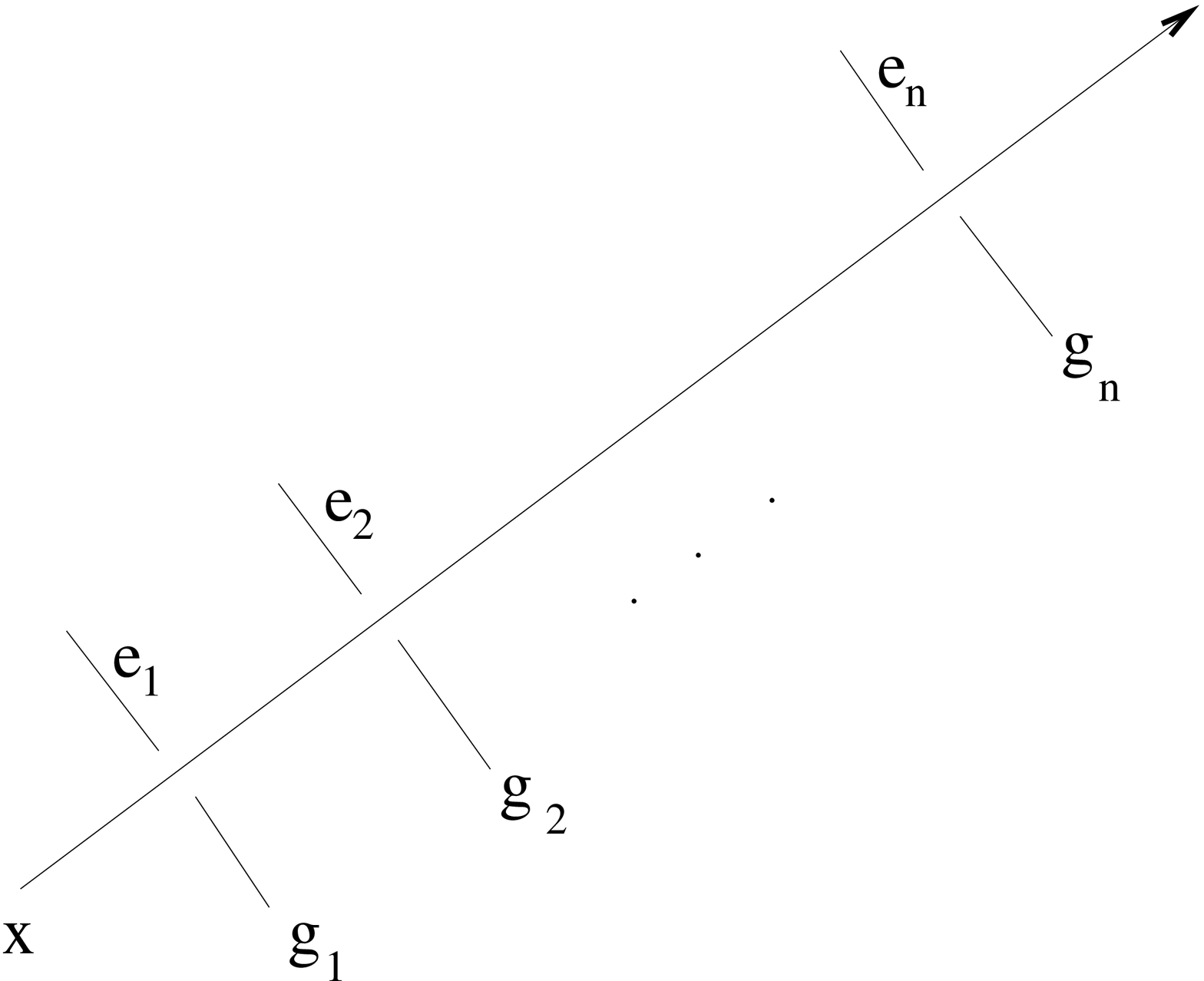,width=1.2in}}
\caption{\label{gauge-plus}}
\end{figure}

\begin{figure}
\centerline{\psfig{figure=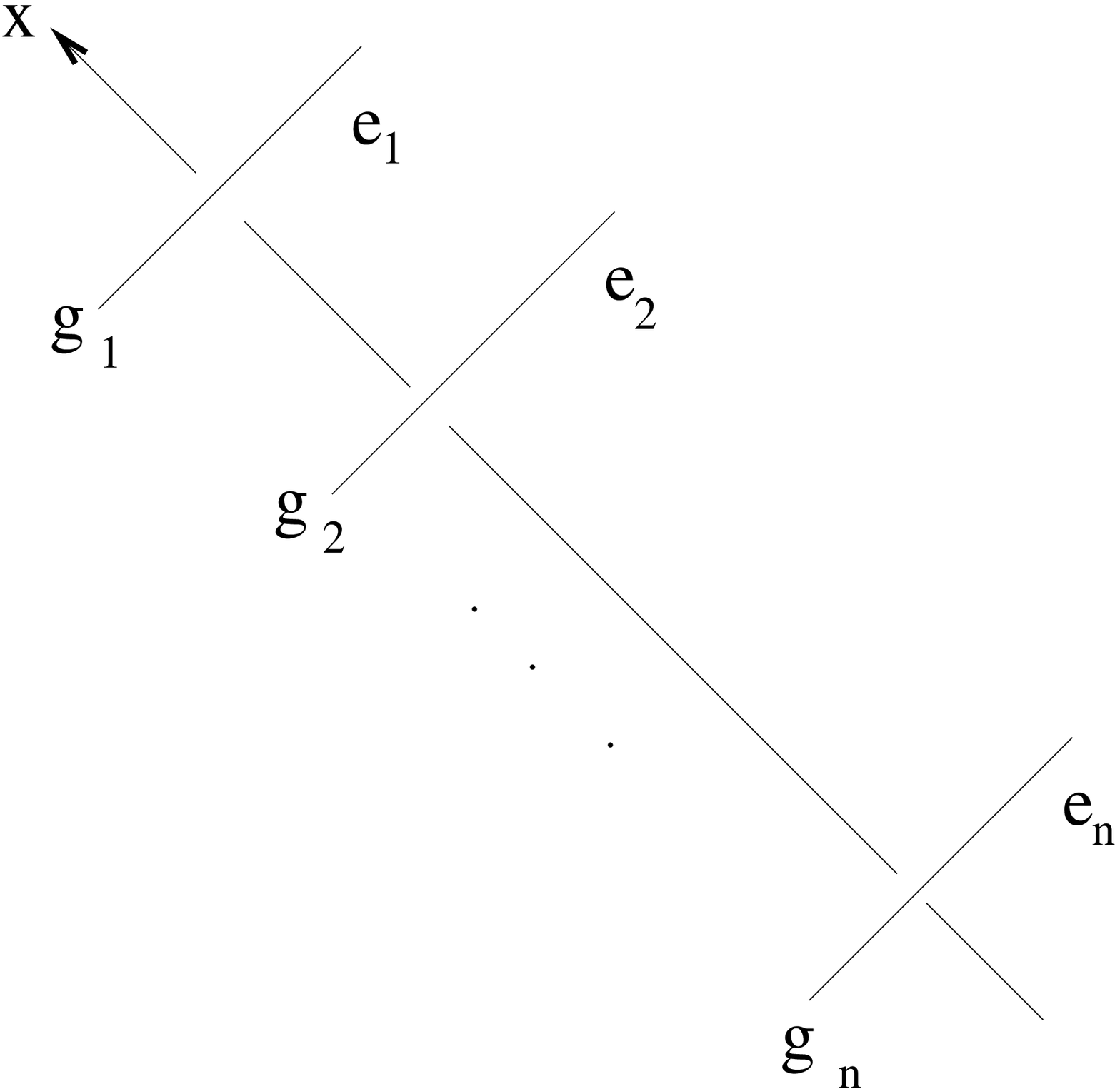,width=1.2in}}
\caption{\label{gauge-minus}}
\end{figure}

\begin{prop} If $g_1,\dots,g_n$ are $G$-colorings of
"lower" edges of diagrams in Figures~\ref{gauge-plus}, \ref{gauge-minus},
then the colorings of "upper" edges are
\[
\ep_i(\tilde{g}_i)=x_+^{(g_{i-1})_-^{-\ep_{i-1}}\dots(g_1)_-^{-\ep_1}}
\ep_i(g_i)\left(x_+^{(g_{i-1})_-^{-\ep_{i-1}}\dots(g_1)_-^{-\ep_1}}\right)^{-1}
\]
for the diagram in Fig~\ref{gauge-plus} and
\[
\ep_i(\tilde{g}_i)=\left(
x_-^{-1}\right)^{(g_1)_+^{\ep_1}\dots(g_{i-1})_+^{\ep_{i-1}}}
\ep_i(g_i)\left(\left(x_-^{-1}\right)^{(g_1)_+^{\ep_1}
\dots(g_{i-1})_+^{\ep_{i-1}}}\right)^{-1}
\]
for the diagram in Fig~\ref{gauge-minus}.
\end{prop}

\begin{proof} Let us consider first the diagram from
Fig.\ref{gauge-plus}. Lemma \ref{+cross} implies the following
relations between $G$-colorings:
\[
\ep_i(\tilde{g}_i)=(x_{i-1})_+\ep_i(g_i)(x_{i-1})_+^{-1},\quad
(x_{i})_+=\left((g_{i})_-^{-\ep_{i}}(x_{i-1})_+^{-1}\right)_+^{-1}
\]
From here we deduce that
\[
(x_i)_+=\left(\left(g_i\right)_-^{-\epsilon_i}\dots
\left(g_1\right)_-^{-\epsilon_1}x_+^{-1}\right)_+^{-1}=
x_+^{\left(g_i\right)_-^{-\epsilon_i}\dots
\left(g_1\right)_-^{-\epsilon_1}}
\]
This, together with the formula (\ref{+dress}) imply the first
statement of the proposition.

Similarly, for the $G$-coloring of the diagram
from Fig. \ref{gauge-minus} we have
\[
\ep_i(\tilde g_i)=(x_{i-1})_-^{-1}\epsilon_i(g_i)(x_{i-1})_-,\quad
(x_i)_-=\left((x_{i-1})_-^{-1}(g_{i})_+^{\epsilon_i}\right)_-
\]
This gives the formula for $(x_i)_-$:
\[
(x_i)_-=\left(x_-^{-1}(g_1)_+^{\ep_1}\dots (g_{i})_+^{\ep_{i}}\right)_-
=\left(\left(x_-^{-1}\right)^{(g_1)_+^{\ep_1}\dots (g_{i})_+^{\ep_{i}}}
\right)^{-1}
\]
Together with the formula (\ref{-dress}) this implies the last
statement of the proposition and completes the proof.

\end{proof}

\subsection{Gauge invariance of the functor $F$}

Let $R^{X,Y}(x,y)$ be a holonomy $R$-matrix acting in $X\otimes Y$
where $X,Y\in \mathfrak X$. Define the operator
\[
\sigma(R^{X, Y}(x,y))=P^{YX}R^{X,Y}(x,y)P^{XY}: \ Y\otimes X\to X\otimes Y
\]

Let $t$ be an $\mathfrak X$-decorated $G$-colored tangle. Assume that
\[
t\in Hom_{{\widetilde{\mathcal FT}}(G,{\mathfrak X})}
((\epsilon_1,Y_1,y_1)\otimes \dots\otimes (\epsilon_n,Y_n,y_n),
(\sigma_1,Z_1,z_1)\otimes\dots\otimes (\sigma_k,Z_k,z_k))
\]
Denote by $F(t)$ the value of functor $F$ described in
Subsection~\ref{invariant}.

\begin{figure}
\centerline{\psfig{figure=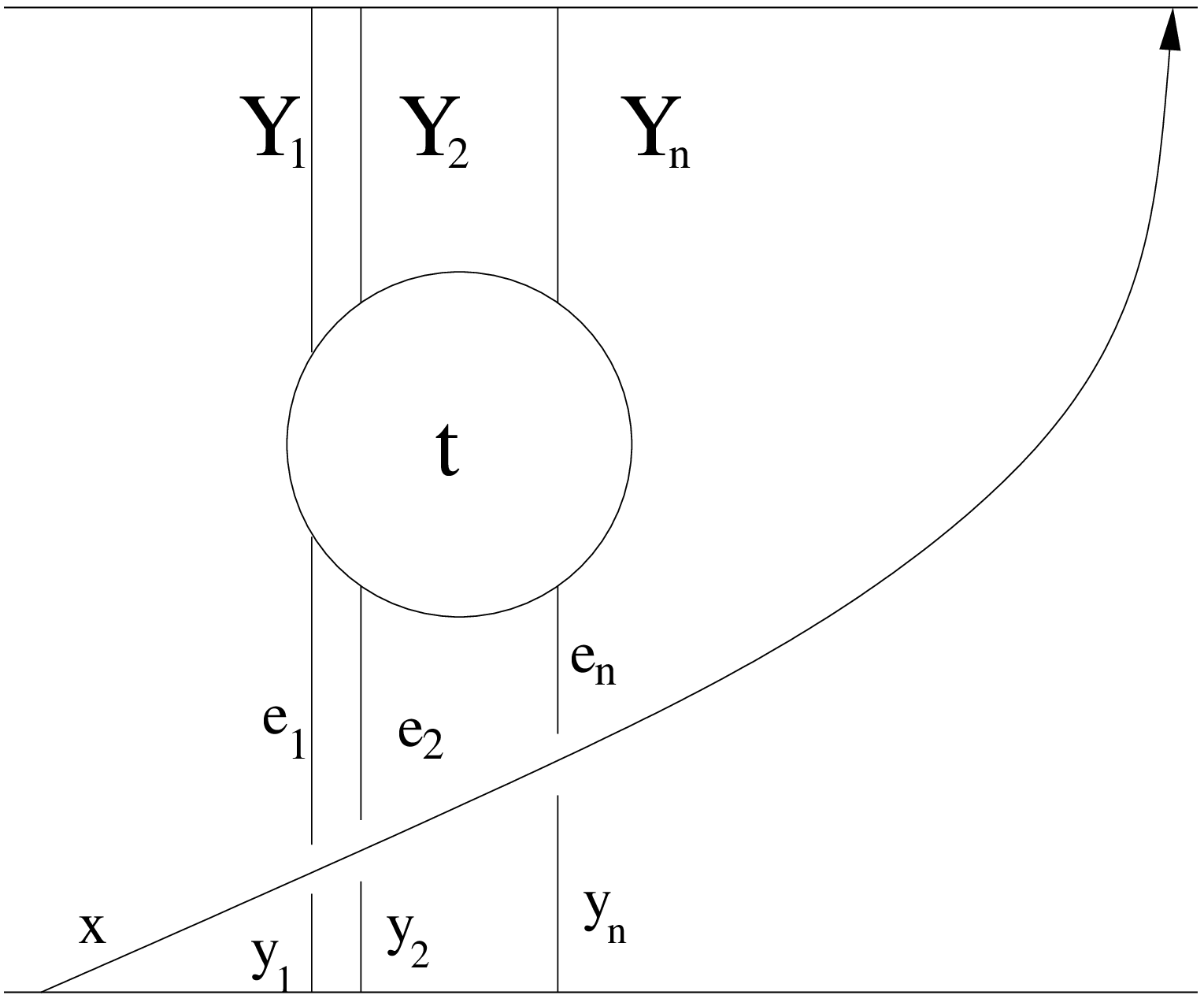,width=1.2in}}
\caption{\label{gauge-invp}}
\end{figure}

\begin{figure}
\centerline{\psfig{figure=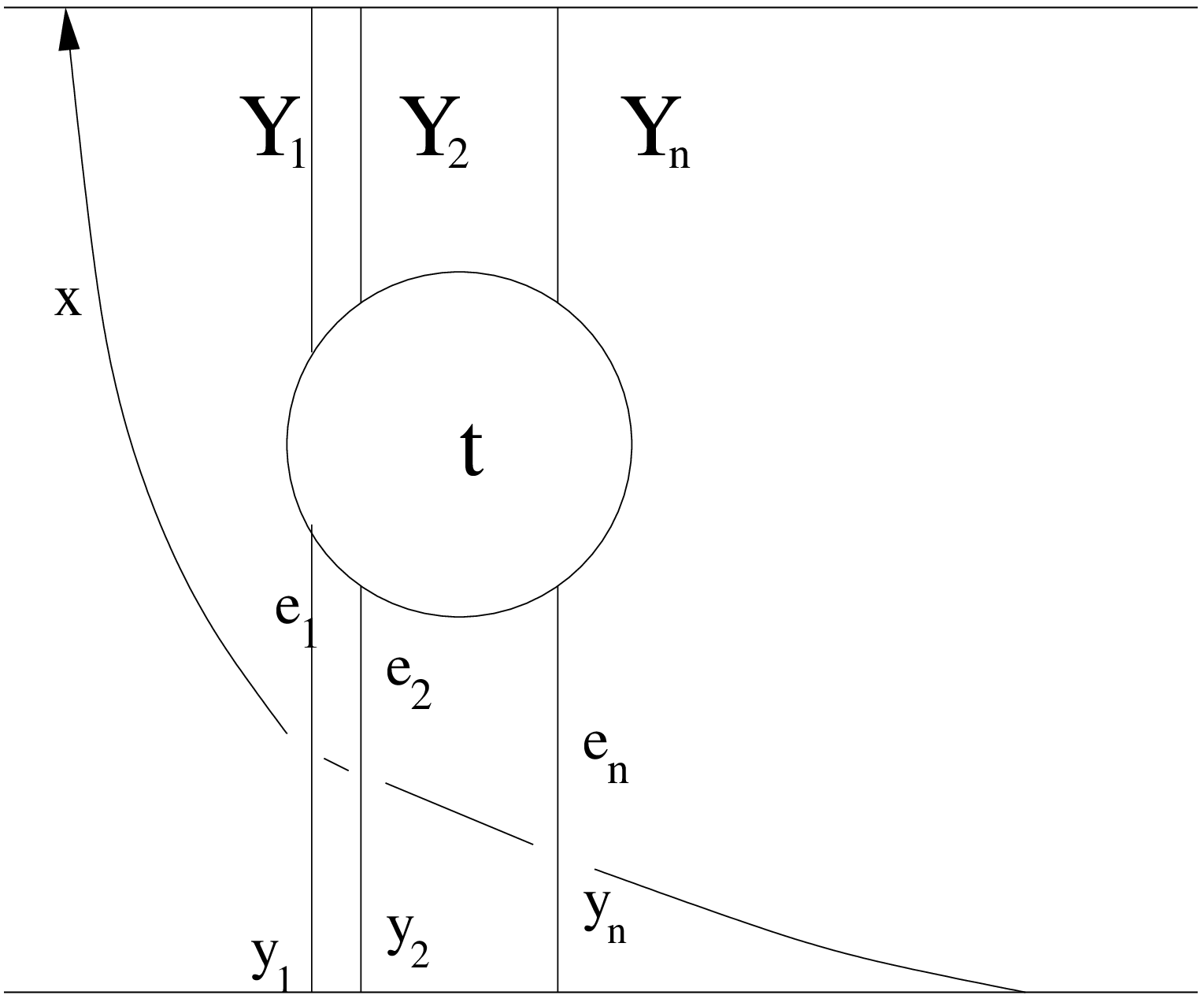,width=1.2in}}
\caption{\label{gauge-invm}}
\end{figure}

\begin{thm} Let $t^x$ be an $\mathfrak X$-decorated $G$-colored tangle
with the $G$-coloring obtained by the gauge action of $x\in G$ on
the $G$-coloring of the the tangle $t$.  Then the following
relations hold for $x_{\pm}\in G_\pm$
\begin{multline}\label{gauge-}
(1\otimes
F(t^{x_+}))R_{0n}^{\e_n,n}(\{y\},x)^{-1}\dots
R_{01}^{\e_1,1}(\{y\},x)^{-1}=\\
R_{0k}^{\sigma_k,k}(\{z\},x)^{-1}\dots
R_{01}^{\sigma_1,1}(\{z\},x)^{-1}(1\otimes  F(t))
\end{multline}
\begin{multline}\label{gauge+}
(F(t^{x_-})\otimes 1)\tilde{R}_{1,n+1}^{\e_1,1}(\{y\},x)^{-1}\dots
\tilde{R}_{n,n+1}^{\e_n,n}(\{y\},x)^{-1}=\\ R_{1,k+1}^{\sigma_1,1}(\{z\},x))^{-1}
\dots R_{k,k+1}^{\sigma_k,k}(\{z\},x)^{-1} (F(t)\otimes 1)
\end{multline}
Both sides in (\ref{gauge-}) are linear maps $X\otimes Y_1\otimes\dots\otimes Y_n\to
X\otimes Z_1\otimes\dots\otimes Z_k$ with
\[
R_{0i}^{(+,i)}(\{y\},x)=\sigma(R^{Y_i,X}(\alpha_{x_-}(y)_i,x_i))_{0i}
\]
\[
R_{0i}^{(-,i)}(\{y\},x)=\sigma((R^{Y_i,X}(\alpha_{x_-}(y)_i,x_{i-1})^{t_1})^{-1})_{0i}
\]
Here we used the notation
\[
(a\otimes b)_{0i}=a_0\otimes 1\otimes \dots\otimes b_i\otimes \dots \otimes 1\in End(
X\otimes Y_1\otimes \dots \otimes Y_n)
\]
Both sides in (\ref{gauge+}) are linear maps $ Y_1\otimes\dots\otimes Y_n\otimes X\to
Z_1\otimes \dots\otimes Z_k\otimes X$ with
\[
\tilde{R}_{i,n+1}^{(+,i)}(\{y\},x)=
\sigma(R^{X,Y_i}(x_{i-1},\alpha_{x_+}(y)_i))_{i,n+1}
\]
\[
\tilde{R}_{i,n+1}^{(-,i)}(\{y\},x)=
\sigma((R^{X,Y_i}(x_{i},\alpha_{x_+}(y)_i)^{-1})^{t_2})_{i,n+1}
\]
where
\[
(a\otimes b)_{i,n+1}=1\otimes\dots\otimes a_i\otimes\dots\otimes 1\otimes
b_{n+1}
\]
The maps $\alpha_{x_\pm}$ are defined in (\ref{alpha+}) (\ref{alpha-}),
and the elements $x_i, \tilde{x_i}, x^i, \tilde{x^i}$ are determined by
the recurrence relations
\[
x_{i-1}=(\alpha_{x_-}(y)_i)_-^{-\epsilon_i}x_i(\alpha_{x_-}(y)_i)_-^{\epsilon_i}
\]
\[
x^{i-1}=(y_i)_+^{\epsilon_i}x_i(y_i)_+^{-\epsilon_i}
\]
with $x_0=x_+$ and $x^0=x_-$. Replacing $y$ by $z$ and $n$ by $k$ we will
have the definition of the remaining $R$-matrices.
\end{thm}
Proof. First let us prove the equation (\ref{gauge-}). Consider
the diagram which is represented on  the Fig. \ref{gauge-invm}.
The value of the functor $F$ on this diagram coincide with the value of
this functor on the diagram which obtained from it by moving
the component colored by $x$ under the tangle to top of the diagram.
This equality is exactly the equality (\ref{gauge-}). Similarly, moving
the component colored by $x$ on the diagram from
Fig. \ref{gauge-invp} we prove the second equality.
Evaluation of the functor $F$ on these diagrams is an easy routine.
Formulae from Appendix \ref{useful} are useful.

Q.E.D.

Thus, the functor $F$ is an invariant of tangles with flat $G$-connection
in the complement which is transforms as it is described above with
respect to the action of the gauge group.

\begin{cor} If $t$ is a link, the invariant $F(t)$ is a scalar,
i.e. it takes values in the base field which was assumed $\Bbb C$.
It is invariant with respect to the gauge group action on the flat
connection and therefore is a function ( or a section of a line
bundle) on the moduli space of flat $G$-connections on the
complement to the tangle $t$.
\end{cor}

\section{Conclusion}
We described a construction of invariants of tangles with flat
connections in the complement which is based on a system of
nondegenerate and cross-nondegenerate holonomy $R$-matrices. In
the second paper of this series we will show how to construct
holonomy $R$-matrices from the representation theory of quantum groups
at roots of 1. We will also show how to use these holonomy
$R$-matrices to construct invariants of 3-manifolds with flat
$G$-connection in it.

\section{Appendix}

\subsection{}
\begin{figure}
\centerline{\psfig{figure=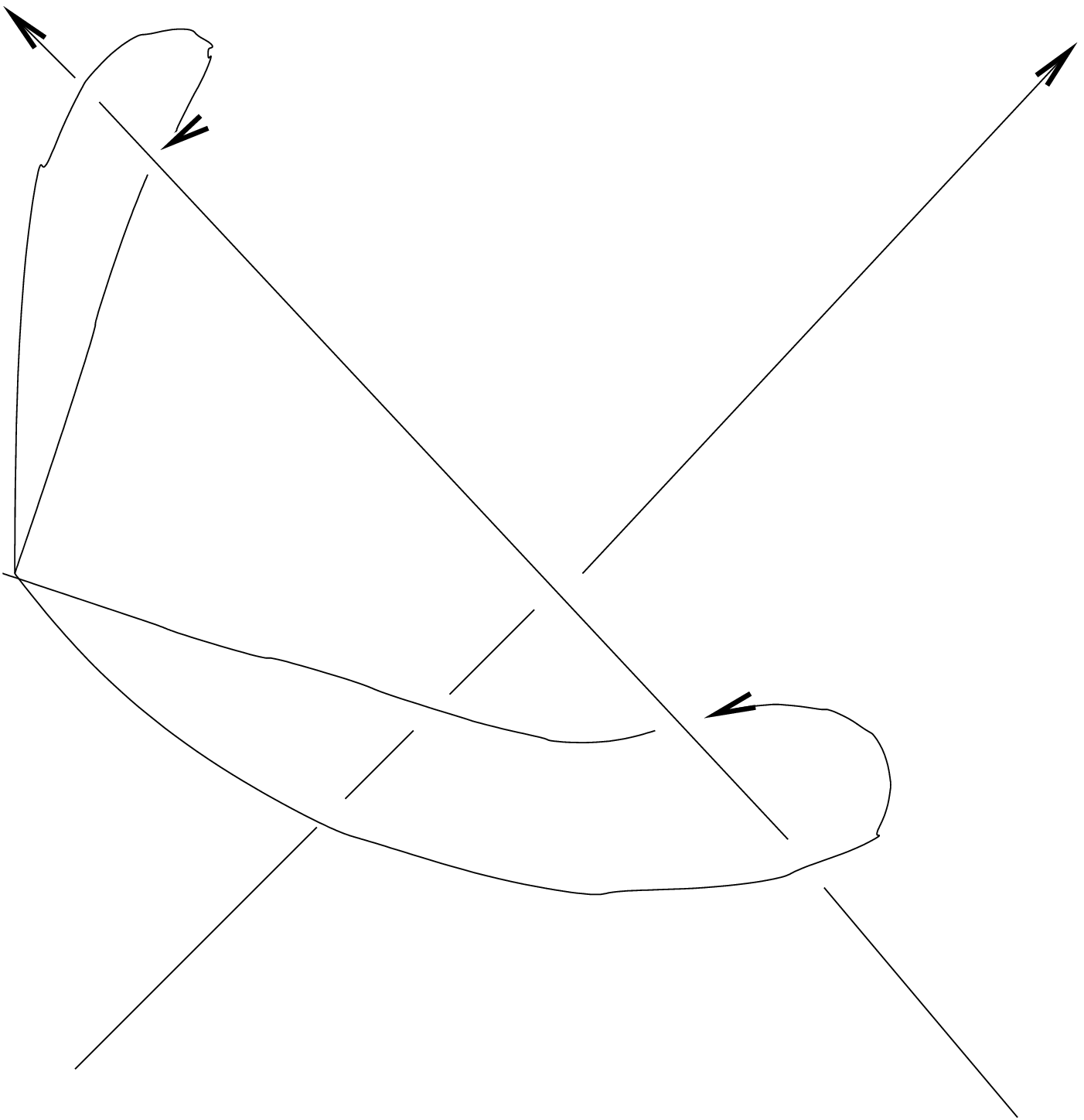,width=1.2in}}
\caption{\label{mon-1}}
\end{figure}

\begin{figure}
\centerline{\psfig{figure=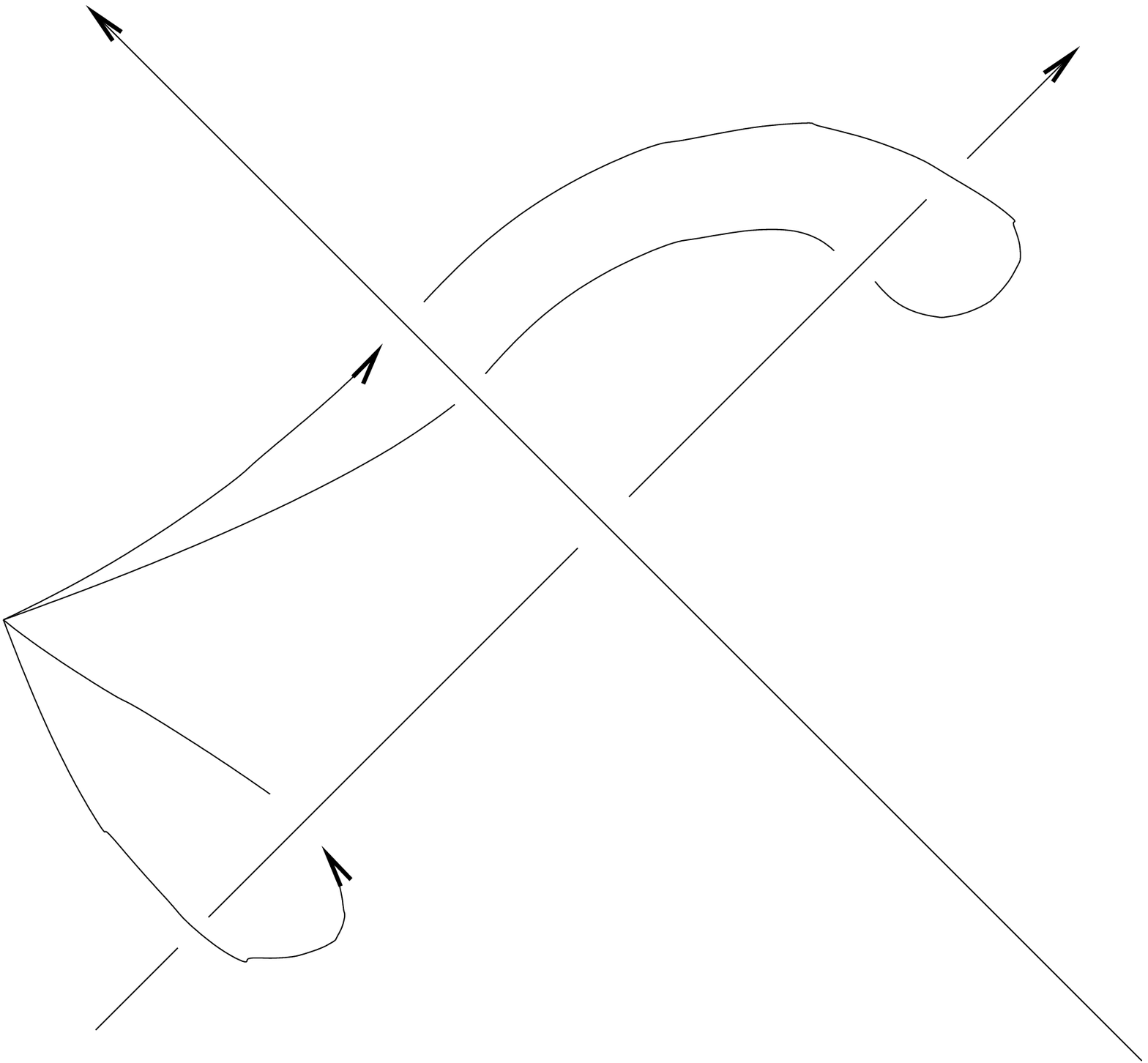,width=1.2in}}
\caption{\label{mon-2}}
\end{figure}

\begin{figure}
\centerline{\psfig{figure=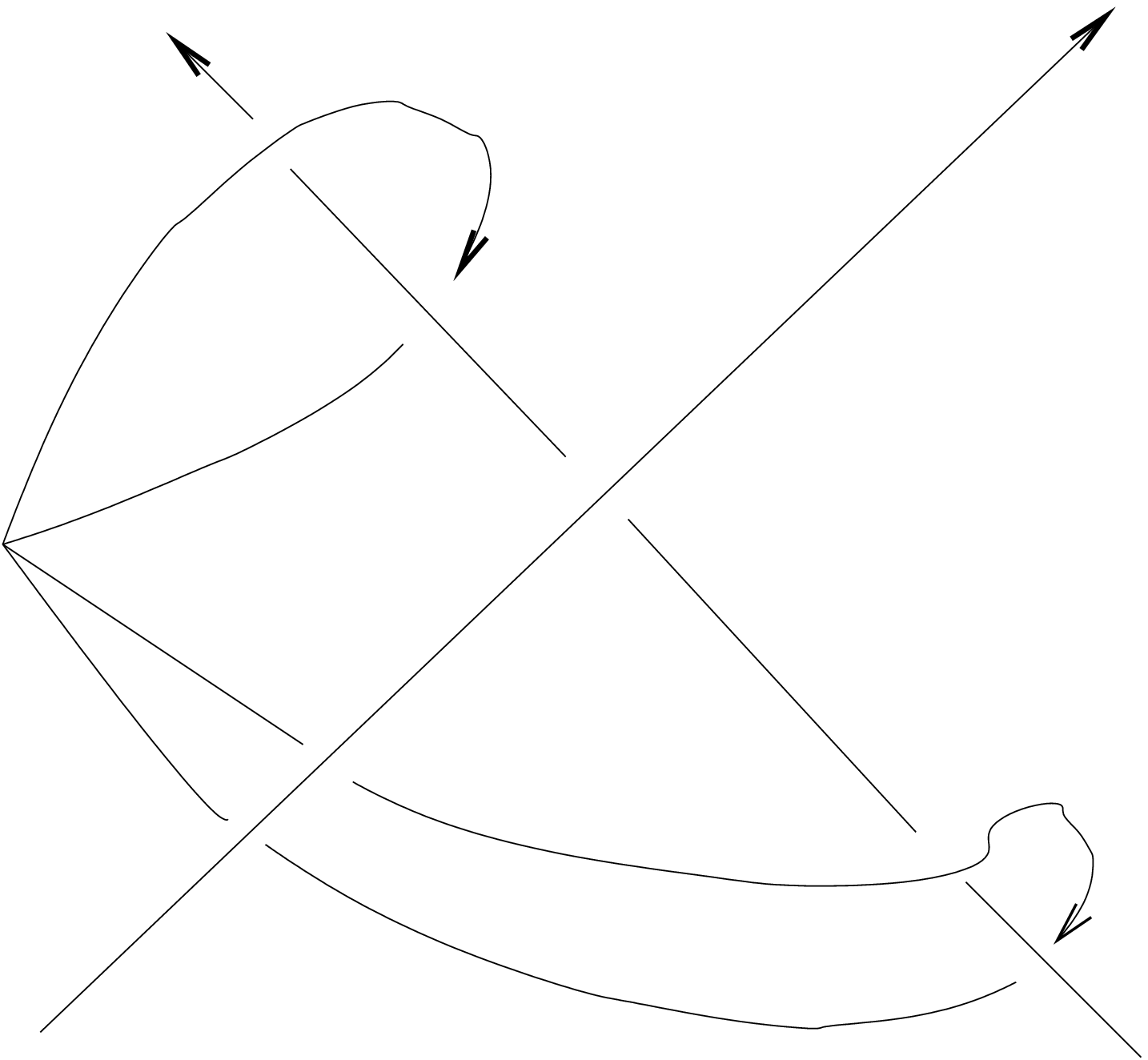,width=1.2in}}
\caption{\label{mon-3}}
\end{figure}

\begin{figure}
\centerline{\psfig{figure=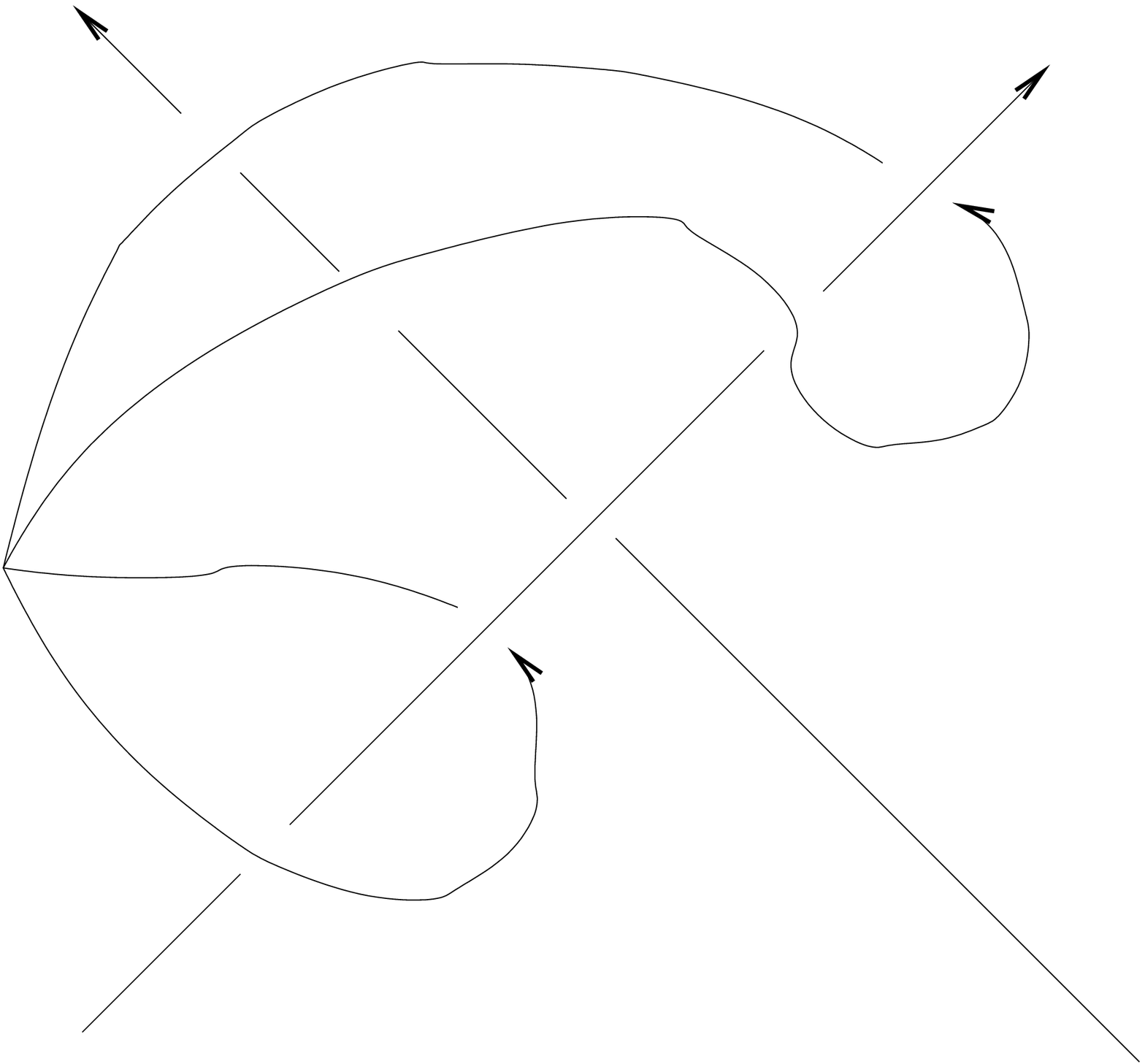,width=1.2in}}
\caption{\label{mon-4}}
\end{figure}

Let us show that the crossing rules from Fig. \ref{crossing-rule} imply
the rules in the definition of the $G$-coloring. Indeed, consider
pairs of pathes from Fig. \ref{mon-1}-\ref{mon-4}. Isotopy equivalence of
these pathes imply the equality of corresponding monodromies. Computing
monodromies along these pathes according to crossing rule from
Fig. \ref{crossing-rule} we obtain the definition of the $G$-coloring.

\subsection{}\label{useful}
\begin{figure}
\centerline{\psfig{figure=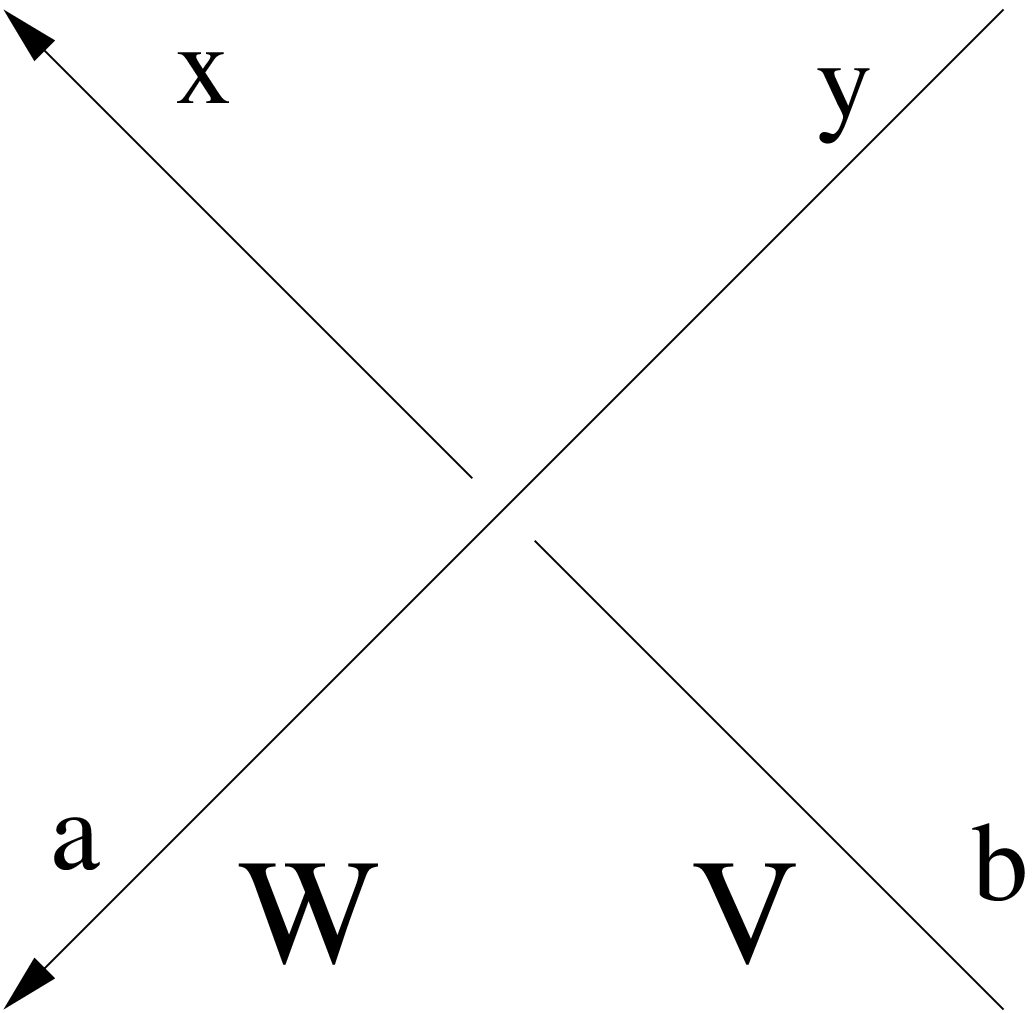,width=1.2in}}
\caption{\label{t-1}}
\end{figure}

\begin{figure}
\centerline{\psfig{figure=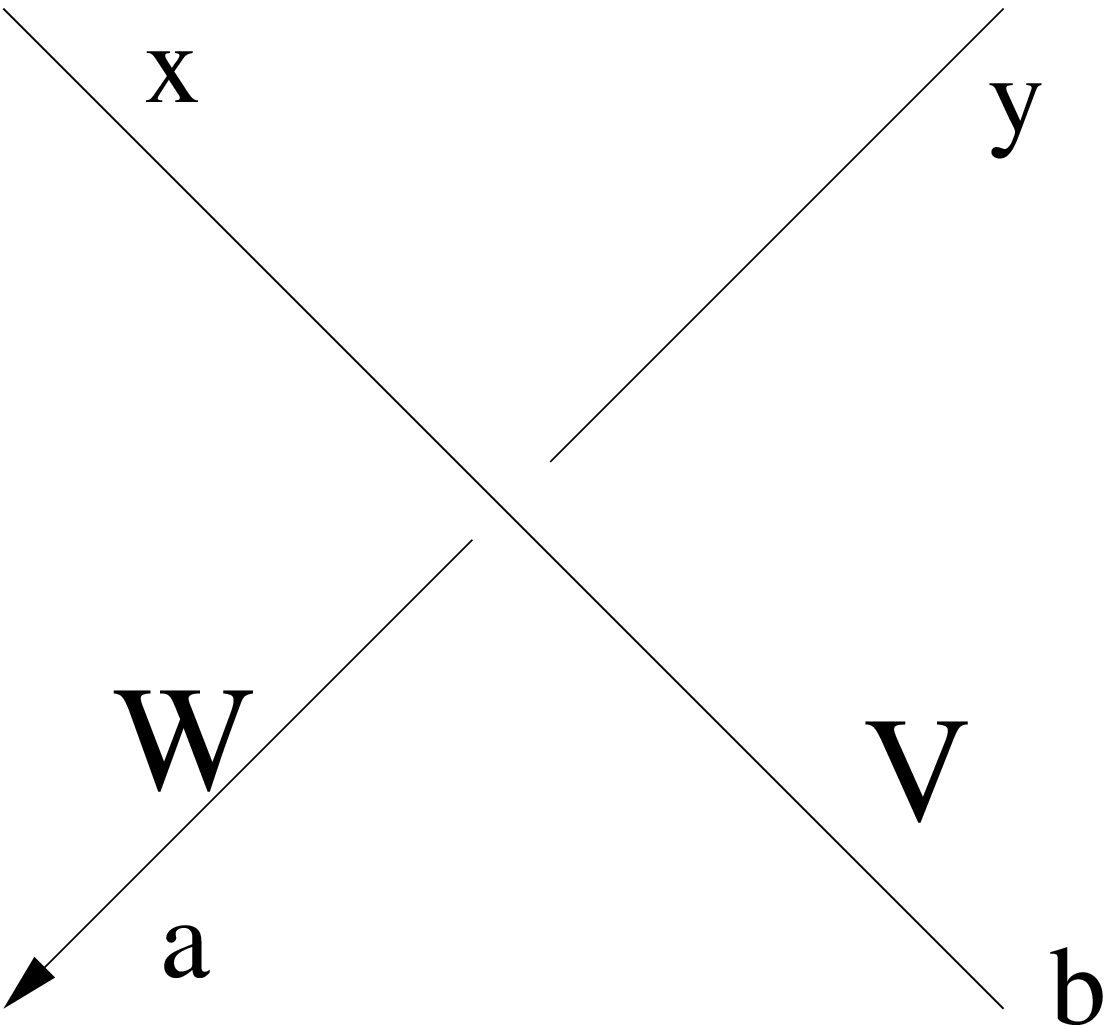,width=1.2in}}
\caption{\label{t-2}}
\end{figure}

\begin{figure}
\centerline{\psfig{figure=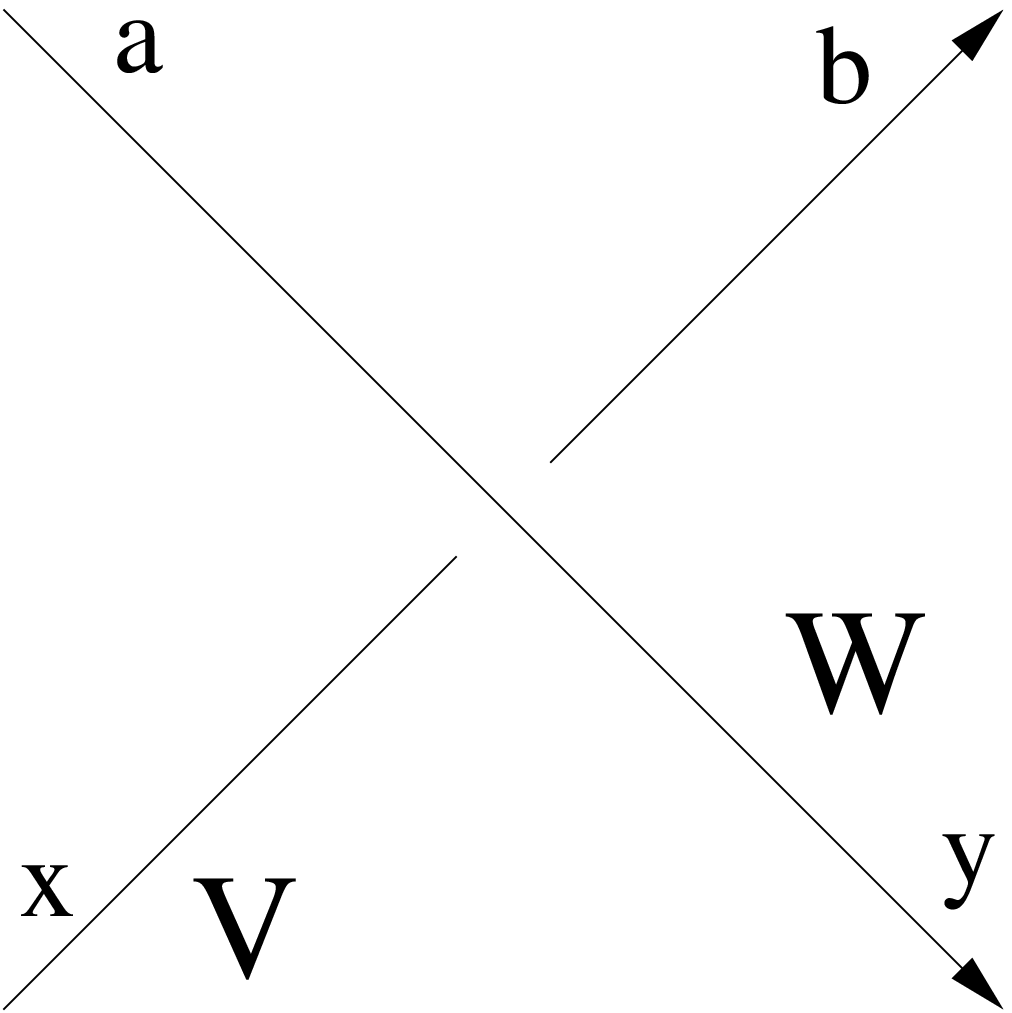,width=1.2in}}
\caption{\label{t-1-inv}}
\end{figure}

\begin{figure}
\centerline{\psfig{figure=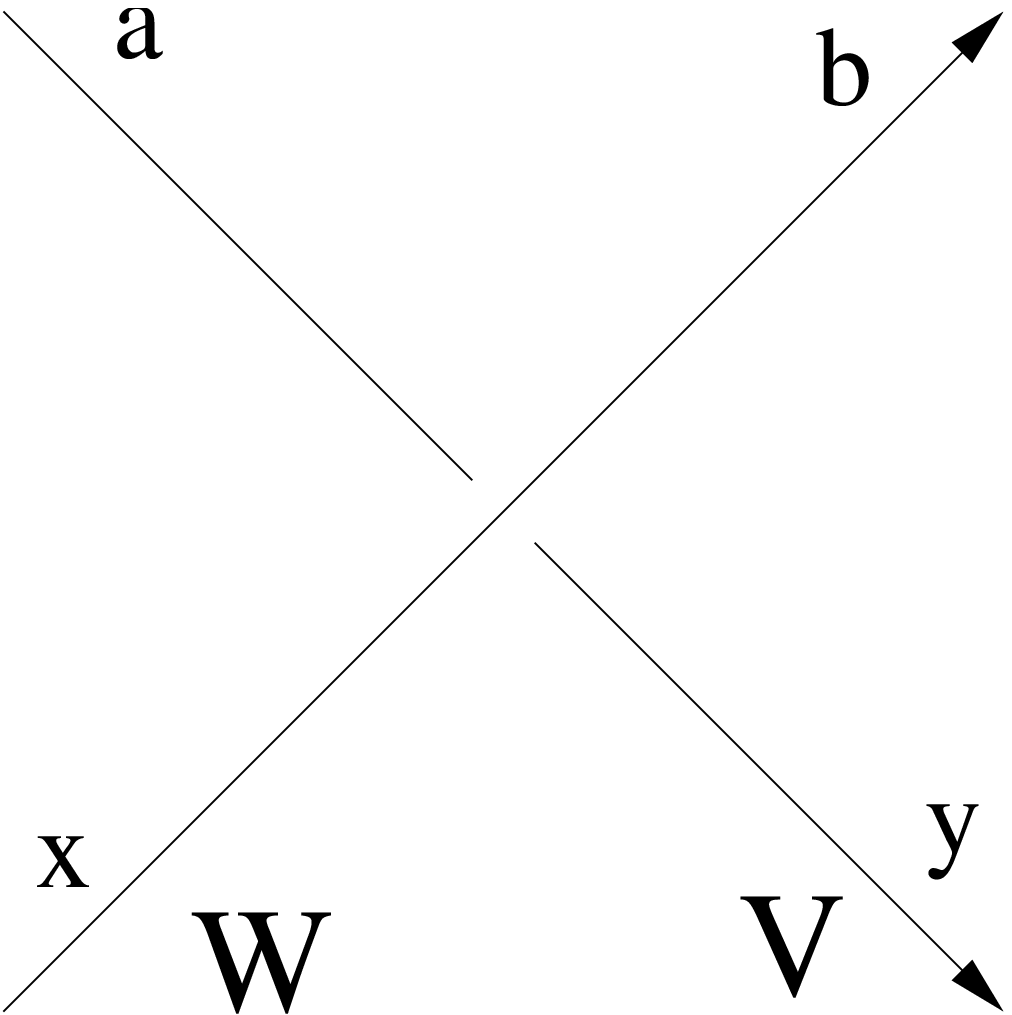,width=1.2in}}
\caption{\label{t-2-inv}}
\end{figure}

Let $t_1, t_2$ and $t_1^{-1}, t_2^{-1}$ be diagrams from Fig.
\ref{t-1}-\ref{t-2-inv}. The functor $F$ can be evaluated on them
using the second framed Redemeister. As a result we have the following
linear maps:
\[
F(t_1)=((R^{VW}(b,y)^{-1})^{t_2})^{-1}P^{W^*V}:W^*\otimes V\to V\otimes
W^*
\]
\[
F(t_2)=P^{W^*V}((R^{WV}(a,x))^{t_1})^{-1}:W^*\otimes V\to V\otimes
W^*
\]
\[
F(t_1^{-1})=F(t_1)^{-1}=P^{VW^*}(R^{VW}(b,y)^{-1})^{t_2}
\]
\[
F(t_2^{-1})=F(t_2)^{-1}=(R^{WV}(a,x))^{t_1}P^{VW^*}
\]

\end{document}